\documentclass[11pt,a4paper]{article}

\usepackage{amsmath,amssymb,amsfonts,amsthm,latexsym,tikz}
\usepackage{graphicx, bm, bbm, yfonts, colortbl}
\usepackage[citecolor=blue,urlcolor=blue,hyperindex,breaklinks]{hyperref}
\usepackage{subcaption,float}
\usepackage{MnSymbol}
\usepackage{verbatim}

\newcommand {\IR}{\mathbb{R}}
\newcommand{\R}{\mathbb{R}}

\newcommand{\re}{\mathbb{R}}

\newcommand {\Z}{{\mathbb Z}}

\newcommand{\cB}{\mathcal B}

\newcommand{\cF}{\mathcal F}

\newcommand{\tf}{\mathcal{F}}

\newcommand{\cI}{\mathcal I}
\newcommand{\cJ}{\mathcal J}

\newcommand{\cK}{\mathcal K}

\newcommand{\cL}{\mathcal L}

\newcommand{\tL}{{\sigma, b}}

\newcommand{\half}{\frac{1}{2}}
\newcommand{\beq}{\begin{equation}}
\newcommand{\eeq}{\end{equation}}
\newcommand{\beqn}{\begin{equation*}}
\newcommand{\eeqn}{\end{equation*}}

\newtheorem{stat}{Statement}[section]
  \newtheorem{prop}[stat]{Proposition}
  
  \newtheorem{thm}[stat]{Theorem}
  \newtheorem{lemma}[stat]{Lemma}
  \newtheorem{remark}[stat]{Remark}

\numberwithin{equation}{section}
\numberwithin{subsection}{section}

\begin{document}


{	\centering 
	
	
	
	
	
	
	{ \scshape \bf \large Asymptotic behavior of the
stochastic heat equation over large intervals}\footnote[1]{Version June 22, 2026


 Asymptotics for boundary value equations.

%
} 
	
	
	
	
	
	
	

	\vspace{0.5 in} 
	
	
			

	{\scshape {\large David Candil, Robert C.~Dalang}\\ \'Ecole Polytechnique F\'ed\'erale de Lausanne} 
	\vspace{0.5\baselineskip}
	
		{\scshape and}
		\vspace{0.5\baselineskip}
			
	{\scshape{\large Marta Sanz-Sol\'e}\\ University of Barcelona} 
	
	\vspace{0.5 in} 

	
	
	
	
	
	



\thispagestyle{empty}

~




}


\begin{abstract}
We consider a nonlinear stochastic heat equation 
on $[0,T]\times [-L,L]$, driven by a space-time white noise $W$, with a given initial condition $u_0: \R \to \R$ and three different types of (vanishing) boundary conditions: Dirichlet, Mixed and Neumann. We prove that as $L\to\infty$, the random field solution  at any space-time position converges in the $L^p(\Omega)$-norm ($p\ge 1$) to the solution of the stochastic heat equation on $\R$ (with the same initial condition $u_0$), and we determine the (near optimal) rate of convergence. The proof relies on estimates of 
differences between the corresponding Green's functions on $[-L, L]$ and the heat kernel on $\R$, and on a space-time version of a Gronwall-type lemma. 
\end{abstract}

\section{Introduction}
\label{s1}

In this paper, we are interested in comparing the random field solutions, as $L \to \infty$, of the stochastic heat equation on $D_L:=[-L, L]$ with appropriate boundary conditions, and the solution on $\R$, when both stochastic partial differential equations (SPDEs) have the same coefficients and the same initial condition. The setting is as follows:
For $\cL = \frac{\partial}{\partial t} - \frac{\partial^{2}}{\partial x^{2}}$, $T>0$, $u_0:\R \to \R$ and $L>0$, we consider the SPDE
\begin{equation}
\label{s1-0}
\begin{cases}
 \cL\, u_L(t, x) = \sigma(t,x,u_L(t, x))\, \dot W(t, x) & \\
 \qquad \qquad \qquad \qquad + b(t, x, u_L(t, x)), &(t, x) \in\, ]0,T] \times D_L,\\
 u_L(0, x) = u_0(x), &x\in D_L,
 \end{cases}
\end{equation}
with one of the following three types of vanishing boundary conditions: for $t > 0$,
\begin{description}
\item[{\rm Dbc}]  (Dirichlet boundary conditions) $u_L(t, -L) = u_L(t, L) = 0$;
  \item[{\rm Mbc}] (Mixed boundary conditions) ~~~$u_L(t, -L) = \frac{\partial}{\partial x} u_L(t, L) = 0$;
 \item[{\rm Nbc}]  (Neumann boundary conditions) $\frac{\partial}{\partial x} u_L(t, -L) = \frac{\partial}{\partial x} u_L(t, L) = 0$.
 \end{description}
 We also consider the stochastic heat equation on $\re$:
 \begin{equation}
\label{s1-1}
\begin{cases}
\cL\, u(t, x) = \sigma(t,x,u(t, x))\, \dot W(t, x) & \\
  \qquad \qquad \qquad \qquad + b(t, x, u(t, x)), &(t, x) \in\, ]0,T] \times \R ,\\
 u(0, x) = u_0(x), &x\in \R . \\
 \end{cases}
\end{equation}
Both equations are driven by the same space-time white noise $W = (W(A),\, A\in\mathcal{B}^f_{\IR_+\times \IR })$, where $\mathcal{B}^f_{\IR_+\times \IR }$ stands for the $\sigma$-field of Borel subsets of $\IR_+\times \IR$ with finite Lebesgue measure. 

Consider the random field solutions $u_L(t,x)$ and $u(t,x)$ to \eqref{s1-0} and \eqref{s1-1}, respectively, for $(t,x)\in[0,T]\times \R$ such that $x \in D_L$. Under suitable conditions which will be specified later, we prove that the $L^p(\Omega)$-norm of the difference $u_L(t,x)-u(t,x)$ goes to zero as $L\to\infty$, and we obtain a nearly optimal estimate on this difference (see Theorem \ref{t-1} below). The estimate is ``nearly optimal" in the sense that in various special cases that include, for instance, the Gaussian case and the case of additive noise ($\sigma \equiv 1$), a slightly sharper estimate can be obtained: see Sections \ref{s1.0} and \ref{rd08_25s1}.


Approximations of deterministic parabolic PDEs on the whole space by PDEs on a bounded domain with Dirichlet boundary conditions appear for instance in \cite{H-R-S-W-2013,L-L-1996}, in the context of  numerical solutions to problems arising in mathematical finance. The error related to these approximations is called the {\em localization error}. 
Numerical approximations are also the motivation for the papers by Gerencsér and Gyöngy \cite{G-G-2017,G-G-2018}, where  localization errors are studied specifically. In a more general setting than ours, and with a method relying on a suitable version of the Feynman-Kac formula for stochastic Dirichlet equations, they find the localization error for approximations in Sobolev norms to a class of SPDEs on the whole space by truncated versions on balls (with Dirichlet boundary conditions). We note that their localization error is similar to that given in Theorem \ref{t-1} below. In \cite{G-G-2017}, the convergence rate of a numerical scheme for the SPDE on the whole space is determined by adding two errors: the localization error and, for the truncated SPDE, the error derived from a finite difference approximation for the spatial variable $x$ and the implicit Euler approximation for the temporal variable $t$.

In this article, we focus on localization errors. Instead of a truncation of \eqref{s1-1}, we fix $L>0$ and consider the boundary value SPDE \eqref{s1-0} subject to the three types of boundary conditions Dbc, Mbc and Nbc. The approximation to the heat equation on $\re$ when $L\to \infty$ is obtained through two steps. First, by a precise analysis of the discrepancy in the $L^p$-norms, $p=1,2$, of integrals of the corresponding Green's functions for \eqref{s1-0} and the heat kernel (the fundamental solution corresponding to \eqref{s1-1}); this is carried out in Section \ref{s3}.  Secondly, we  use an explicit computation of a geometric space-time convolution series related to the heat kernel that appears in \cite{chendalang2015-2}, together with a space-time version of a Gronwall-type lemma (see Lemma \ref{a1-l3}).

A similar approach can be used to study SPDEs on $\re^k$, $k\ge1$, driven by a Gaussian noise that is white in time and with a spatially homogeneous covariance \cite[Chapter 4]{candil2022}.
Other possible extensions would include SPDEs driven by strictly parabolic (regular) partial differential operators, and  SPDEs with rough initial conditions (such as unbounded functions or measures): see \cite{chendalang2015-2}, \cite[Section 6.6]{d-ss-2024} and \cite{CCL} for details.

The organization of this paper is as follows.  Section \ref{s1.0} is devoted to the study of simplified linear versions of \eqref{s1-0} and \eqref{s1-1}, and serves to introduce the subject.  In Section \ref{s1.1}, we fix the setting of the paper and state the main theorem (Theorem \ref{t-1}) on upper bounds on the localization error in $L^p(\Omega)$-norm ($p \geq 1$), which is proved in Section \ref{s2}. In Section \ref{rd08_25s1}, we discuss simplified nonlinear versions of \eqref{s1-0} and \eqref{s1-1}, in order to get further insight into the near optimality of Theorem \ref{t-1}, and we also establish a lower bound on the localization error in $L^2(\Omega)$-norm in the cases of linear multiplicative noise and of uniformly positive $\vert \sigma \vert$. In Section \ref{s3}, we establish bounds on the discrepancy, in terms of $L^p$-norms, $p=1, 2$, between the Green's function $\Gamma_L$, corresponding to the stochastic heat equation with boundary conditions, and the heat kernel $\Gamma$. The results are fundamental for the proof of Theorem \ref{t-1}, but they are also of independent interest. The last section is an appendix that contains auxiliary results used throughout the paper, and, in particular, a space-time Gronwall-type lemma inspired by the calculations in \cite{chendalang2015-2} and \cite{candil2022}.


\section{The deterministic and the linear cases}
\label{s1.0}

Before entering into the core of the article, we discuss some special cases of \eqref{s1-0} and \eqref{s1-1}, namely the homogeneous PDE ($\sigma \equiv 0$ and $b \equiv 0$), and the linear SPDE ($\sigma \equiv 1$ and $b \equiv 0$). These two special cases suggest the types of bounds that we will aim for in Section \ref{s1.1} for  the solution of \eqref{s1-1}.


We denote by $\Gamma(t, x)$ the fundamental solution associated to the heat operator on $\R$:
\beq
\label{s1-3}
\Gamma(t,x) = \frac{1}{\sqrt{4\pi t}} \exp\left(-\frac{x^2}{4t}\right) 1_{]0,\infty[}(t),\quad x\in\re,
\eeq
and by $\Gamma_L(t;x,y)$ the Green's function associated to the heat equation \eqref{s1-0} on $D_L$ with one of the three types of boundary conditions Dbc, Mbc and Nbc described in Section \ref{s1}. Whenever we need to differentiate between the three boundary conditions, we write $\Gamma_L^D(t;x,y)$, $\Gamma_L^M(t;x,y)$ and $\Gamma_L^N(t;x,y)$, respectively, instead of $\Gamma_L(t;x,y)$ (see \eqref{2}, \eqref{4} and \eqref{3} for their respective expressions).

 For $L > 0$, $t\in\, ]0,T]$ and $x,y\in D_L$,  let 
\beqn
H_L(t; x, y) = \Gamma(t, x-y) - \Gamma_L(t; x, y),
\eeqn 
and define $K_{t, x, L}$ by 
\beq\label{rd07_18e2}
    K_{t, x, L} = \min\left(\tfrac 1 2, \sqrt{\tfrac{t}{\pi}}\, \tfrac{1}{L-\vert x \vert} \right).
\eeq
Clearly, $K_{t, \pm L, L} = \frac12$. The function $x\mapsto \Phi(x)$ will denote the standard Normal probability distribution function and $\vartheta(a) := \sum_{m=0}^\infty e^{-am^2}$, $a>0$,  will denote the  {\em Theta function} \cite[Chapter 5, Section 3.1]{stein-shak-2003}.

For the two cases Dbc and Mbc (vanishing Dirichlet and Mixed boundary conditions), let $L_0 = 0$, and for the case Nbc (vanishing Neumann boundary conditions), fix $L_0 > 0$.  

In the two next Propositions, we assume that the initial condition $u_0:\re\to \re$ is a bounded Borel function, and we set $\Vert u_0\Vert_\infty = \sup_{x\in \re} |u_0(x)|$. Moreover, unless explicitly stated otherwise, the results are valid for all three of the boundary conditions Dbc, Mbc and Nbc.


 \begin{prop} \label{rd05_04e1}
Fix $T > 0$ and suppose that $\sigma \equiv 0$ and $b\equiv 0$.  

(a) For all $L > L_0$ and $(t, x) \in [0, T] \times D_L$,
 \begin{align*}
   \vert u(t, x) - u_L(t, x))\vert \leq K_{t, x, L} \, \Vert u_0 \Vert_{\infty} \left[\exp\left(-\tfrac{(L-x)^2}{4t} \right) + \exp\left(-\tfrac{(L+x)^2}{4t} \right) \right].
 \end{align*}
 
 (b)  Consider the case of vanishing Dirichlet boundary conditions Dbc. Assume that $m:= \inf_{y \in \R} u_0(y) \geq 0$. Then for all $L > 0$ and all $(t, x) \in [0, T] \times D_L$, 
 \begin{align*}
&\vert u(t, x) - u_L(t, x)\vert \\
 &\qquad \geq m\, \sqrt{\tfrac{t}{\pi}}\, \left[ \tfrac{L-x}{2t+(L-x)^2} \exp\left( -\tfrac{(L-x)^2}{4t}\right) + \tfrac{L+x}{2t+(L+x)^2} \exp\left( -\tfrac{(L+x)^2}{4t}\right)\right]. 
 \end{align*}
\end{prop}

\begin{proof}
(a) Under the assumptions of the proposition, 
\begin{align*}
    u(t, x) = \int_\R dy\, \Gamma(t, x-y)\, u_0(y) \quad\text{and}\quad u_L(t, x) = \int_{D_L} dy\, \Gamma_L(t; x, y)\, u_0(y),
\end{align*}
therefore,
\begin{align*}
    u(t, x) - u_L(t, x)  &=   \int_{D_L^c} dy\, \Gamma(t, x-y)\, u_0(y)  \\
   &\qquad\qquad +  \int_{D_L} dy\, H_L(t;x, y)\, u_0(y) .
\end{align*}

By Lemma \ref{a1-l2}, and Lemmas \ref{s3-l1}, \ref{s3-l3} and \ref{s3-l2} (specifically \eqref{cd-1}, \eqref{cn-10} and \eqref{cn-1}), for each of the three types of boundary conditions Nbc, Mbc and Dbc,
\begin{align*}
   \vert u(t, x) - u_L(t, x) \vert &\leq K_{t, x, L} \, \Vert u_0\Vert_{\infty} \left[\exp\left(-\tfrac{(L-x)^2}{4t} \right) + \exp\left(-\tfrac{(L+x)^2}{4t} \right) \right].
\end{align*}
This proves (a).

(b) In the case of vanishing Dirichlet boundary conditions Dbc, $H_L > 0$ by Lemma \ref{s3-l1}, therefore, 
\begin{align*}
   \vert u(t, x) - u_L(t, x) \vert &\geq  m \int_{D_L^c} dy\, \Gamma(t, x-y) \\
   &= m \int_{\frac{L-x}{\sqrt{2t}}}^{+\infty} dz\, \varphi(z) + \int_{\frac{L+x}{\sqrt{2t}}}^{+\infty} dz\, \varphi(z),
\end{align*}
  where $\varphi(u) := (2 \pi)^{-\half}\exp\left(-\frac{u^2}{2}\right)$ and we have used the change of variables $u = (y-x)/\sqrt{2t}$. 

   We now use the classical Mill's inequality \cite[7.8]{OLBC}
\begin{align}\label{rd05_05e1}
    1 - \Phi(u) \geq \frac{u}{1+ u^2}\, \varphi(u),\qquad u \geq 0,
\end{align}
to see, after simplification, that $ \vert u(t, x) - u_L(t, x) \vert$ is bounded below by
\begin{align*}
 m \left[\sqrt{2t}\, \tfrac{L-x}{2t+(L-x)^2} \, \varphi\left(\tfrac{L-x}{\sqrt{2t}} \right) +  \sqrt{2t}\, \tfrac{L+x}{2t+(L+x)^2} \, \varphi\left(\tfrac{L+x}{\sqrt{2t}} \right)\right].
\end{align*}
This establishes (b).
\end{proof}
\medskip

 \begin{prop} \label{rd04_29e1}
 Fix $T > 0$ and suppose that $\sigma \equiv 1$, $b\equiv 0$.
 
   (a) For all $L > L_0$ and $(t, x) \in [0, T] \times D_L$, 
\begin{align}\notag
    &\Vert u(t, x) - u_L(t, x)\Vert_{L^2(\Omega)}\\
    &\qquad\qquad\leq \left( K_{t, x, L}^2 \Vert u_0\Vert_{\infty}^2 +  \tfrac{K_{t/2, x, L}}{\sqrt{8 \pi t}} +\sqrt{\tfrac{t}{\pi}}\, 2\vartheta\left(\tfrac{4L_0^2}{t}\right) K_{t, x, L} \right)^\half \notag \\
    &\qquad\qquad\qquad \times \left[\exp\left(-\tfrac{(L-x)^2}{4t} \right) + \exp\left(-\tfrac{(L+x)^2}{4t} \right) \right],
 \label{rd05_05e3}
 \end{align}
 with $2\vartheta\left(\tfrac{4L_0^2}{t}\right)$ replaced by $1$ in cases Dbc and Mbc.
 
   (b) For $t > 0$ and $x \in \R$ fixed, there is $L_1 = L_1(t, x) > 0$ such that for all  $L > \max(L_1, \vert x \vert)$,  \begin{align}
 \label{mss05_22-2}
    \Vert u(t, x) - u_L(t, x) \Vert_{L^2(\Omega)} \geq \tfrac{\sqrt{t}}{\sqrt{2 \pi}}\,  \left(1 + \tfrac{L-\vert x\vert}{\sqrt{t}} \right)^{-3/2} \exp\left(- \tfrac{(L-\vert x\vert)^2}{4t} \right) .
 \end{align}
 
   (c) Consider the case of vanishing Dirichlet boundary conditions Dbc. Suppose that $m:= \inf_{y\in \re}u_0(y) \geq 0$. Then for $t > 0$ and $x \in \R$ fixed, there is $L_1 = L_1(t, x) > 0$ such that for all  $L > \max(L_1, \vert x \vert)$,  \begin{align}
 \label{mss05_22-3}
   & \Vert u(t, x) - u_L(t, x) \Vert_{L^2(\Omega)}\notag\\
   &\qquad\quad \ge  
 \left[ \left(m\,  \sqrt{\tfrac{t}{\pi}}\, \tfrac{L-|x|}{2t+(L-|x|)^2}\right)^2 + \sqrt{\tfrac{t}{2\pi}}\,  \left(1 + \tfrac{L-\vert x\vert}{\sqrt{t}} \right)^{-3} \right]^{\frac12} 
  \exp\left(- \tfrac{(L-\vert x\vert)^2}{4t} \right) .
 \end{align}
 \end{prop}

 \begin{proof}
 Under the assumptions of this proposition,
 \begin{align*}
 u(t, x) &= \int_\re dy\ \Gamma(t,x-y) u_0(y) + \int_0^t \int_{\R} \Gamma(t-s, x-y)\, W(ds, dy), \\
  u_L(t, x) &= \int_{D_L} dy\ \Gamma_L(t;x,y) u_0(y) + \int_0^t \int_{D_L} \Gamma_L(t-s; x, y)\, W(ds, dy),
 \end{align*}
 therefore
 \begin{align}
  u(t, x) - u_L(t, x)&=
   \int_{D_L^c} dy\  \Gamma(t,x-y) u_0(y)+ \int_{D_L} dy\ H_L(t;x,y) u_0(y) \notag\\
    &\qquad\qquad + \int_0^t \int_{D_L^c} \Gamma(t-s, x-y)\, W(ds, dy)  \notag\\
     &\qquad\qquad
     + \int_0^t \int_{D_L} H_L(t-s; x,y)\, W(ds, dy).
 \label{rd05_05e2}
 \end{align}
 In particular,
  \begin{align}
  \label{mss-5_22e2}
  &E[(u(t, x) - u_L(t, x))^2]\notag\\
   &\qquad \qquad = \left(\int_{D_L^c} dy\  \Gamma(t,x-y) u_0(y)+ \int_{D_L} dy\ H_L(t;x,y) u_0(y)\right)^2 \notag\\
  &\qquad \qquad \qquad +  \int_0^t ds \int_{D_L^c} dy\,  \Gamma^2(t-s, x-y)
  +  \int_0^t ds \int_{D_L} dy\, H_L^2(t-s; x,y).
 \end{align}
By Proposition 2.1 (a), and by Lemma \ref{a1-l2}, and by Lemmas \ref{s3-l1}, \ref{s3-l3} and \ref{s3-l2} (specifically \eqref{cd-3}, \eqref{cn-12} and \eqref{cn-3}), for each of the three types of boundary conditions Nbc, Mbc and Dbc, 
  \begin{align*}
  &E[(u(t, x) - u_L(t, x))^2] \\
  &\qquad \leq K_{t, x, L}^2\,  \Vert u_0\Vert_{\infty}^2\left[\exp\left(-\tfrac{(L-x)^2}{4t} \right) + \exp\left(-\tfrac{(L+x)^2}{4t} \right) \right]^2 \\
    &\qquad \qquad  + \tfrac{K_{t/2, x, L}}{\sqrt{8 \pi t}}   \left[\exp\left(-\tfrac{(L-x)^2}{2t}\right) + \exp\left(-\tfrac{(L+x)^2}{2t}\right)\right]\\
    &\qquad\qquad +  \sqrt{\tfrac{t}{\pi}}\, 2  \vartheta\left(\tfrac{4L_0^2}{t}\right)K_{t, x, L} \left[\exp\left(-\tfrac{(L-x)^2}{4t}\right) + \exp\left(-\tfrac{(L+x)^2}{4t}\right)\right]^2,
   \end{align*}
 with $2 \vartheta\left(\tfrac{4L_0^2}{t}\right)$ replaced by $1$ in cases Dbc and Mbc. This proves \eqref{rd05_05e3}.

 (b) From \eqref{mss-5_22e2}, for $L > x \geq 0$, we see that
  \begin{align*}
  &E[(u(t, x) - u_L(t, x))^2] \\
    &\qquad \geq \int_0^t ds \int_{D_L^c} dy\, \Gamma^2(t-s, x-y) 
   \geq  \int_0^t ds \int_{L}^{+\infty} dy\, \Gamma^2(t-s, x-y) \\
  &\qquad =  \int_0^t dr \int_{L}^{+\infty} dy\, \Gamma^2(r, x-y) = \int_0^t \tfrac{dr}{\sqrt{8 \pi r}} \int_{\frac{L-x}{\sqrt{r}}}^{+\infty} du\, \varphi(u), 
  \end{align*}
  where $\varphi(u) := (2 \pi)^{-\half}\exp\left(-\frac{u^2}{2}\right)$ and we have used the change of variables $u = (y-x)/\sqrt{r}$, and so
\begin{align*}
  E[(u(t, x) - u_L(t, x))^2] \geq \tfrac{1}{2 \sqrt{2\pi}} \int_0^t \tfrac{dr}{\sqrt{r}}\, \left(1 - \Phi\left( \tfrac{L-x}{\sqrt{r}}\right) \right).
\end{align*}

   We use again Mill's inequality \eqref{rd05_05e1}
to see, after simplification, that $E[(u(t, x) - u_L(t, x))^2]$ is bounded below by
\begin{align*}
   \tfrac{1}{2\sqrt{2\pi}} \int_0^t  dr\, \tfrac{L-x}{r + (L-x)^2} \, \varphi\left( \tfrac{L-x}{\sqrt{r}}\right)
   \geq \tfrac{1}{4\pi}\, \tfrac{L-x}{t + (L-x)^2}\, \int_0^t dr\, \exp\left(- \tfrac{(L-x)^2}{2r} \right).
\end{align*}
We now apply Lemma \ref{rd04_29l1} in the Appendix to complete the proof of\eqref{mss05_22-2}. 

   (c) Finally, \eqref{mss05_22-3} is a consequence of \eqref{mss-5_22e2}, Proposition \ref{rd05_04e1} (b) and the calculations that led to \eqref{mss05_22-2}.
 \end{proof}

\section{Hypotheses, preliminaries and main result}
\label{s1.1}

In this section, we suppose that the initial condition $x \mapsto u_0(x)$ is Borel and bounded and we also assume the following properties ($\mathbf H_{\sigma, b}$) 
on the coefficients $\sigma$ and $b$:
\medskip

\noindent ($\mathbf H_{\sigma, b}$): 
\smallskip

\noindent{\em Measurability and adaptedness.} The functions $\sigma$ and $b$  are defined on $[0,T]\times\re\times \IR\times \Omega$ with values in $\IR$ and are jointly measurable, that is, $\mathcal{B}_{[0,T]}\times \mathcal{B}_\re \times\mathcal{B}_{\IR}\times \tf$-measurable. These two functions are also adapted to $(\tf_s,\, s\in[0,T])$ (the completed natural filtration associated with $W$), that is, for fixed $s \in [0,T]$, $(y,z,\omega)\mapsto \sigma(s,y,z;\omega)$ and $(y,z,\omega)\mapsto b(s,y,z;\omega)$ are $\cB_\re\times \cB_{\IR}\times \cF_s$-measurable.
\medskip

\noindent {\em Global Lipschitz condition.}
There exists $\bar C:=\bar C_T\in \IR_+$ such that for all $(s,y,\omega)\in [0,T]\times \re\times \Omega$ and $z_1,z_2\in \IR$,
\beqn
\vert \sigma(s,y,z_1;\omega)-\sigma(s,y,z_2;\omega)\vert + \vert b(s,y,z_1;\omega)-b(s,y,z_2;\omega)\vert  \le \bar C\, |z_1-z_2|.
\eeqn

\noindent{\em Uniform linear growth.} There exists a constant $\bar c :=\bar c_T\in \IR_+$ such that for all $(s,y,\omega)\in [0,T]\times \re \times \Omega$ and all $z\in \IR$,
\beqn
\vert \sigma(s,y,z;\omega)\vert + \vert b(s,y,z;\omega)\vert \le \bar c \left(1+|z|\right).
\eeqn
\medskip

In the sequel, without loss of generality, we will assume for simplicity that $\bar C=\bar c$ and we denote this constant $C_\tL$.


The random field solutions to \eqref{s1-0} and \eqref{s1-1} are the processes $(u_L(t,x),\, (t,x)\in[0,T]\times D_L)$ and $(u(t,x),\ (t,x)\in[0,T]\times \re)$, satisfying, for $(t, x) \in [0, T] \times D_L$,
\begin{align}
\label{ul}
u_L(t,x) = I_{0,L}(t,x) &+ \int_0^t \int_{D_L} \Gamma_L(t-s; x,y)\, \sigma(s,y,u_L(s,y))\, W(ds,dy)\notag\\
&+  \int_0^t ds \int_{D_L} dy\, \Gamma_L(t-s; x,y)\, b(s,y,u_L(s,y)),
\end{align}
and for $(t, x) \in [0, T] \times \R$,
\begin{align}
\label{u}
    u(t, x) =  I_0(t,x) &+ \int_0^t \int_\R \Gamma(t-s, x-y)\, \sigma(s,y,u(s,y))\, W(ds,dy)\notag\\
&+  \int_0^t  ds \int_\R dy\, \Gamma(t-s, x-y)\, b(s,y,u(s,y)),
\end{align}
where
\beq
\label{initial-cond}
I_{0,L}(t,x) = \int_{D_L}\Gamma_L(t;x,y) u_0(y)\,dy ,\quad  I_0(t,x)=\int_\re\Gamma(t,x-y) u_0(y)\,dy.
\eeq

The stochastic integrals in \eqref{ul} and \eqref{u} are defined in the sense of Walsh \cite{walsh}; see also \cite[Chapter 2]{d-ss-2024} for the properties of this integral.
The inequalities listed in Lemma \ref{l1}, along with the assumptions $(\mathbf H_{\sigma, b})$ and those on the initial condition $u_0$, yield 
the existence and uniqueness of the random field solutions $(u_L(t,x))$ and $(u(t,x))$ in \eqref{ul} and \eqref{u}, and moreover, for all $p>0$, these processes satisfy
\beq
\label{s1-4}
   \sup_{(t,x) \in [0,T]\times D_L} E[\vert u_L(t, x)\vert^p] < \infty,\quad   \sup_{(t,x) \in [0,T]\times \re}  E[\vert u(t, x)\vert^p] < \infty
\eeq
 (see e.g.~\cite[Theorem 4.2.1]{d-ss-2024}). 

The main result of this paper is the following:

\begin{thm}
\label{t-1}
Fix $T > 0$ and $p \geq 1$. Let $u_0: \R \to \R$ be Borel and bounded and set $\Vert u_0\Vert_\infty= \sup_{x \in \re}|u_0(x)|$.
Assume that the functions $\sigma$ and $b$ satisfy $({\mathbf H_{\sigma, b}})$. Let $(u_L(t, x),\, (t, x) \in [0,T] \times D_L)$ be the random field solution to the SPDE \eqref{s1-0}and let  $(u(t, x),\, (t, x) \in [0,T] \times \R)$ be the random field solution to \eqref{s1-1}.

 For the two cases Dbc and Mbc (vanishing Dirichlet and Mixed boundary conditions), let $L_0 = 0$, and for the case Nbc (vanishing Neumann boundary conditions), fix $L_0 > 0$.  
Then there is a constant $c = c(p, T, C_\tL, L_0) < \infty$
such that, for all $L > L_0$ and $(t, x) \in [0, T] \times D_L$,
\begin{align}\nonumber
  \Vert u(t, x) - u_L(t, x) \Vert_{L^p(\Omega)} &\leq c\, (1 + \Vert u_0 \Vert_\infty) \\
  &\qquad \qquad \times \left[\exp\left(-\tfrac{(L-x)^2}{8t} \right) + \exp\left(-\tfrac{(L+x)^2}{8t} \right) \right].
\label{s1-5}
\end{align}
In particular, for all $(t, x) \in [0, T] \times \R$,
\beq
\label{s1-6}
   \lim_{L \to \infty} u_L(t, x) =  u(t, x) \quad\text{ in }  L^p(\Omega).
\eeq
\end{thm}

\begin{remark}
\label{rem-to-rd08_09t1}

(a) The near optimality of the bound in \eqref{s1-5} is illustrated by the examples of Section \ref{s1.0}. This question will also be discussed further in the examples of Section \ref{rd08_25s1}.

(b) In the case Nbc, it is not possible to take $L_0 = 0$. Indeed, we notice that $0 \in D_L$ for all $L > 0$, and in the case $\sigma\equiv 1$ and $b \equiv 0$, by the semigroup property \eqref{l1-semigroup} and the first formula in \eqref{3},
\begin{align*}
   E[u_L^2(t, 0)] &\geq \int_0^t ds \int_{-L}^L dy \, \left(\Gamma_L^N(t-s; 0, y)\right)^2 = \int_0^t ds \, \Gamma_L^N(2(t-s); 0, 0) \\
    &\geq  \int_0^t ds \, \tfrac{1}{2L}
    = \tfrac{t}{2L},
\end{align*}
therefore
\begin{align*}
    \lim_{L \downarrow 0} E[u_L^2(t, 0)] =\infty.
\end{align*}

(c) If $u_0: \re\times \Omega\rightarrow \re$ is a measurable random variable that is independent of the noise $\dot W$ and such that $\sup_{x\in\re}\Vert u_0(x)\Vert_{L^p(\Omega)}<\infty$, then the conclusion of
Theorem \ref{t-1} still holds, with $\Vert u_0 \Vert_{\infty}$ there replaced by $\sup_{x\in\re}\Vert u_0(x)\Vert_{L^p(\Omega)}$.
\end{remark}

\section{Proof of Theorem \ref{t-1}}
\label{s2}

We begin this section by proving a strengthening of the first statement in \eqref{s1-4}, namely that the $L^p$-moments of the random field solutions to \eqref{ul} are not only bounded uniformly in $(t,x)\in[0,T]\times D_L$, but also in $L$. 

Recall that the initial condition $u_0$ is a bounded Borel function; hence, from \eqref{initial-cond}, we have
\beq
\label{bI0}
\Vert I_0\Vert_\infty:=\sup_{(t,x) \in [0,T]\times \re} |I_0(t,x)| = \Vert u_0\Vert_\infty<\infty
\eeq
and
\beq
\label{bI0-bis}
\sup_{L >  0}\Vert I_{0,L}\Vert_\infty:=\sup_{L > 0}\, \sup_{(t,x) \in [0,T] \times D_L}\, |I_{0,L}(t,x)| \le \Vert u_0 \Vert_\infty <\infty.
\eeq
Indeed, $\int_{\re}\Gamma(t,x-y)\,dy = 1$ and this implies \eqref{bI0}.
By Lemma \ref{l1}, in the three instances of boundary conditions Dbc, Mbc and Nbc, we have
\begin{align*}
\Vert I_{0,L}\Vert_\infty &\le \sup_{(t,x) \in [0,T]\times D_L} \int_{D_L} \Gamma_L(t;x,y)\, |u_0(y)|\,dy\\
& \le \Vert u_0 \Vert_\infty  \sup_{(t,x) \in [0,T]\times D_L} \int_{D_L} \Gamma_L(t;x,y)\,dy \le \Vert u_0 \Vert_\infty
\end{align*}
(see \eqref{l1-a4},  \eqref{l1-c4} and \eqref{l1-b0-bis}). 

\begin{lemma}
\label{rd08_09l1}
The assumptions are the same as in Theorem \ref{t-1}. 
For all $p>0$, there is a constant $C(p, T) < \infty$ such that
\beq
\label{s1-7}
    \sup_{L > L_0}\, \sup_{(t,x) \in [0, T]\times D_L} \Vert u_L(t, x)\Vert_{L^p(\Omega)} \leq C(p, T)\, (1 + \Vert u_0 \Vert_\infty) < \infty.
\eeq
\end{lemma}

\begin{proof}
It suffices to prove the lemma for $p \geq 2$, so we fix $p \geq 2$ and use \eqref{ul} to see that for $x \in D_L$,
\begin{align*}
   E\left[\vert u_L(t, x)\vert^p \right]&\leq 3^{p-1}\left(\vert I_{0,L}(t,x)\vert^p
   + E\left[\vert \cI_L(t, x)\vert^p\right] + E\left[\vert \cJ_L(t, x)\vert^p\right]\right),
\end{align*}
where
\begin{align*}
     \cI_L(t, x) &= \int_0^t \int_{D_L} \Gamma_L(t-s; x,y)\, \sigma(s,y,u_L(s,y))\, W(ds,dy), \\
     \cJ_L(t, x)&= \int_0^t ds \int_{D_L} dy\, \Gamma_L(t-s; x,y)\, b(s,y,u_L(s,y)).
\end{align*}
Using Burkholder's inequality (see e.g.~\cite[Equation (2.2.14)]{d-ss-2024})
and then Hölder's inequality for $\cI_L$, and Hölder's inequality for $\cJ_L$, we obtain
\begin{align*}
   &E\left[\vert u_L(t, x)\vert^p\right] \\
    &\qquad\leq C_p\, \left(\vphantom{ \left[\int_0^t ds \, J_{1,L}(s) \right]^{\frac{p}{2} - 1}  } \Vert u_0\Vert_\infty^p \right.\\
    &\qquad\qquad +\left[\int_0^t ds \, J_{1,L}(s) \right]^{\frac{p}{2} - 1} \int_0^t ds \, J_{1,L}(t-s)
        \left(1 + \sup_{x \in D_L} E\left[\vert u_L(s, x)\vert^p\right]\right) \\
    &\qquad\qquad+ \left. \left[\int_0^t ds \, J_{2,L}(s) \right]^{p - 1} \int_0^t ds \, J_{2,L}(t-s) 
        \left(1 + \sup_{x \in D_L} E\left[\vert u_L(s, x)\vert^p\right]\right)\right),
\end{align*}
where $C_p$ depends on $p$ and $C_\tL$ but not on $L$, and
\beqn
J_{1,L}(s) = \sup_{x\in D_L} \int_{D_L} \Gamma_L^2(s; x,y)\, dy,\qquad
J_{2,L}(s) = \sup_{x\in D_L} \int_{D_L} \Gamma_L(s; x,y)\, dy.
\eeqn

We now consider separately each of the three types of boundary conditions.
\medskip

\noindent{\em Case Dbc.}\ Here $\Gamma_L(t;x,y)=\Gamma_L^D(t;x,y)$, given in \eqref{2}. From \eqref{l1-a2} and \eqref{l1-a4}, we see that for any $s\in[0,T]$,
\beq\label{bound-J}
J_{1,L}(s) + J_{2,L}(s) \leq \tfrac{1}{\sqrt{8 \pi s}}+1 =: J(s),
\eeq
and from \eqref{l1-a3} and  \eqref{l1-a4}, we deduce that
\beqn
\int_0^t ds\, J_{1,L}(s) \le \left(\tfrac{t}{2\pi}\right)^\half \qquad\text{and}\qquad \int_0^t ds\,J_{2,L}\le t.
\eeqn
Let $g_L(t) = \sup_{x\in D_L}  E[\vert u_L(t, x)\vert^p]$. These considerations show that there is $\tilde C(p, T) < \infty$ such that
\beqn
  g_L(t) \leq \tilde C(p, T) \left(\Vert u_0 \Vert_\infty^p + \int_0^T ds \, (1 + g_L(s))\, J(t-s) \right).
\eeqn
By the Gronwall-type Lemma \ref{a1-l3-usual}, we conclude that there is $C(T, p) < \infty$ such that for all $L > L_0 = 0$ and all $t \in [0, T]$, $g_L(t) \leq C(p, T)(1+ \Vert u_0\Vert_\infty^p)$. Since this right-hand side does not depend on $t$ or $L$, \eqref{s1-7} is proved.
\medskip

\noindent{\em Case Mbc.}\ Here $\Gamma_L(t;x,y)=\Gamma_L^M(t;x,y)$, given in \eqref{4}. 
Using \eqref{l1-c2} and \eqref{l1-c4}, we see that for any $s\in[0,T]$,
\beqn
J_{1,L}(s) + J_{2,L}(s) \leq \tfrac{1}{\sqrt{2\pi s}}+1 =:  J(s).
\eeqn
Moreover, \eqref{l1-c3} and \eqref{l1-c4} imply that
\beqn
\int_0^t ds\,J_{1,L}(s) \le \left(\tfrac{2t}{\pi}\right)^\half \qquad\text{and}\qquad \int_0^t J_{2,L}(s) < t.
\eeqn
The remainder of the proof of \eqref{s1-7} is similar to the previous case.
\medskip

\noindent{\em Case Nbc.}\ Here $\Gamma_L(t;x,y)=\Gamma_L^N(t;x,y)$, given in \eqref{3}. 
Assume that $L\ge L_0 > 0$.
Using  \eqref{l1-b5-bis} and \eqref{l1-b0-bis}, for any $s\in[0,T]$, we see that
\beqn
J_{1,L}(s) + J_{2,L}(s) \leq   \tfrac{1}{2L} + \tfrac{1}{\sqrt{2 \pi s}} + 1 \leq 1 + \tfrac{1}{2 L_0} + \tfrac{1}{\sqrt{2 \pi s}}  =: J(s).
\eeqn 
Further, \eqref{l1-b0-bis} implies that
\beqn
    \int_0^t ds\,J_{2,L}(s)\le t \leq T,
\eeqn
and by \eqref{int_J_2-tri}, 
\beqn
    \int_0^t ds\,J_{1,L}(s)\le\tfrac{t}{2 L}+ \sqrt{\tfrac{2t}{\pi}} \leq  \tfrac{T}{2 L_0} + \sqrt{\tfrac{2T}{\pi}}.
\eeqn
Then, we argue as in the two previous cases (except that here, $L_0 > 0$) to complete the proof of \eqref{s1-7}.
%
\end{proof}
\smallskip

\noindent{\em Proof of Theorem \ref{t-1}.}\
It suffices to prove \eqref{s1-5} for $p \geq 2$, so we fix such a $p$ throughout the proof. Using \eqref{ul} and \eqref{u},  we write
\beq
\label{s1-9}
    u(t, x) - u_L(t, x) = \sum_{i=0}^3 A_i(t, x),
\eeq
where
\begin{align*}
    A_0(t, x) &= \int_{D_L^c} dy\,  \Gamma(t,x-y)\, u_0(y), \\
    A_1(t, x) &= \int_{D_L} dy\, H_L(t; x, y)\, u_0(y) \\
     &\qquad + \int_0^t \int_{D_L} H_L(t-s; x, y)\\
     &\qquad\qquad\qquad \times    [\sigma(s,y,u_L(s,y))\, W(ds,dy) 
         + b(s,y,u_L(s,y))\, ds dy\, ],\\
          A_2(t, x) &=  \int_0^t \int_{D_L^c} \Gamma(t-s, x-y)\, [\sigma(s,y,u(s,y))\, W(ds,dy)  \\
      &\qquad\qquad\qquad\qquad  + b(s,y,u(s,y))\, ds dy\, ], \\
    A_3(t, x) &= \int_0^t \int_{D_L} \Gamma(t-s, x-y)\\
     &\qquad\qquad\qquad \times    [(\sigma(s,y,u(s,y)) - \sigma(s,y,u_L(s,y))) \, W(ds,dy) \\
     &\qquad\qquad\qquad\qquad  + (b(s,y,u(s,y)) - b(s,y,u_L(s,y)))\, ds dy\, ].
\end{align*}

Set
\beq
\label{aela}
    a_L(t, x) =  \exp\left(-\tfrac{(L-x)^2}{8t} \right) + \exp\left(-\tfrac{(L+x)^2}{8t} \right)
\eeq
and
\beqn
    f_L(t, x) = \Vert u(t, x) - u_L(t, x) \Vert_{L^p(\Omega)}.
\eeqn

By the triangle inequality,
\beq
\label{s1-10}
    f_L(t, x) \leq\sum_{i=0}^3 \Vert A_i(t, x) \Vert_{L^p(\Omega)}.
\eeq
First, we will prove that
\beq
\label{zero-one-two}
\sum_{i=0}^2 \Vert A_i(t, x) \Vert_{L^p(\Omega)} \le  C_1\,  (1+ \Vert u_0 \Vert_\infty)\, \sqrt{K_{t, x, L}} \,  a_L(t/2,x).
\eeq
In the cases Dbc and Mbc, $ C_1 = C_1(p, T,C_\tL)$, while in the case Nbc, $ C_1= C_1(p, T, C_\tL, L_0)$. We will use several times the fact that $K_{t, x, L} \leq \sqrt{K_{t, x, L}}$ (because $K_{t, x, L} \leq 1$), and that $a_L(t/2, x) \leq a_L(t, x)$.

In order to prove \eqref{zero-one-two}, we apply \eqref{a1-l2-1} to obtain
\beq
\label{bound0}
  \Vert A_0(t, x) \Vert_{L^p(\Omega)} =\vert A_0(t, x)\vert \leq  \Vert u_0\Vert_\infty\, K_{t, x, L} \, a_L(t/2, x).
\eeq
Let
\begin{align*}
\Vert A_{1,1}(t, x) \Vert_{L^p(\Omega)} &= \left\vert \int_{D_L}dy\, H_L(t;x,y) u_0(y)\right\vert,\\
\Vert A_{1,2}(t, x) \Vert_{L^p(\Omega)} &=\left\Vert  \int_0^t \int_{D_L} H_L(t-s; x, y) \sigma(s,y,u_L(s,y))\, W(ds,dy) \right\Vert_{L^p(\Omega)},\\
\Vert A_{1,3}(t, x) \Vert_{L^p(\Omega)} &=\left\Vert  \int_0^t ds \int_{D_L}dy\, H_L(t-s; x, y) b(s,y,u_L(s,y))\right\Vert_{L^p(\Omega)}.
\end{align*}

Using Lemmas \ref{s3-l1}--\ref{s3-l2} (see \eqref{cd-1}, \eqref{cn-10} and \eqref{cn-1}), in each of the cases Dbc, Mbc and Nbc, we have 
\beqn
\Vert A_{1,1}(t, x) \Vert_{L^p(\Omega)} \le \Vert u_0\Vert_{\infty} \int_{D_L} dy\, \vert H_L(t;x,y)\vert \le \Vert u_0\Vert_{\infty} \, K_{t, x, L}\, a_L(t/2,x).
\eeqn

Applying Burkholder's inequality and using \eqref{s1-7}, we see that
\begin{align*}
\Vert A_{1,2}(t, x) \Vert_{L^p(\Omega)}    & \le c_1(p,t,C_\tL)\, (1+ \Vert u_0 \Vert_\infty)\, \left[\int_0^t ds \int_{D_L} dy\, H_L^2(t-s;x,y)\right]^{\half} \\
    &\le c_2 \, (1+ \Vert u_0 \Vert_\infty) \, \sqrt{K_{t, x, L}}\ a_L(t/2,x)
\end{align*}
in each of the cases Dbc, Mbc and Nbc, because of the estimates \eqref{cd-3}, \eqref{cn-12} and \eqref{cn-3}  of Lemmas \ref{s3-l1}--\ref{s3-l2}, respectively. In the cases Dbc and Mbc, $c_2:=c_2(p,T,C_\tL)$, while for Nbc, $c_2:=c_2(p,T,C_\tL, L_0)$ ($L_0>0$).

Applying Minkowski's inequality and using again \eqref{s1-7}, we deduce that
\begin{align*}
\Vert A_{1,3}(t, x) \Vert_{L^p(\Omega)}  &\le \tilde C(p, T, C_\tL)\, (1+ \Vert u_0 \Vert_\infty) \int_0^t ds \int_{D_L} dy\, \vert H_L(t-s;x,y)\vert \\
     &\le c_3(p, T, C_\tL)\,   (1+ \Vert u_0 \Vert_\infty)\, K_{t, x, L}\, a_L(t/2,x),
\end{align*}
 because of the estimates \eqref{cd-2}, \eqref{cn-11} and \eqref{cn-2}.

By the triangle inequality, 
\begin{align}\nonumber
 \Vert A_1(t, x) \Vert_{L^p(\Omega)} &\leq \sum_{j=1}^3 \Vert A_{1,j}(t, x) \Vert_{L^p(\Omega)}  \\
  & \leq  C_1\, (1+ \Vert u_0 \Vert_\infty)\, \sqrt{K_{t, x, L}} \ a_L(t/2, x).
 \label{bound1}
\end{align}

In the cases Dbc and Mbc, $c_2 := c_2(p,T,C_\tL)$ and $C_1 := C_1(p,T,C_\tL)$, while for Nbc, $c_2 := c_2(p,T,C_\tL, L_0)$ and $C_1 := C_1(p,T,C_\tL, L_0)$ ($L_0 > 0$).

With similar arguments, but using \eqref{s1-4} for $u(t,x)$ (instead of \eqref{s1-7}), and the inequalities \eqref{a1-l2-2} and \eqref{a1-l2-1}, we see that
\begin{align}
    & \Vert A_2(t, x) \Vert_{L^p(\Omega)}\notag\\       
    &\qquad \leq  \bar C(p, T, C_\tL)\, (1+ \Vert u_0 \Vert_\infty) \Bigg(\left[\int_0^t ds \int_{D_L^c} dy\, \Gamma^2(t-s, x-y)\right]^{\frac{1}{2}} \notag \\
    &\qquad\qquad\qquad\qquad\qquad\qquad\qquad\quad + \int_0^t ds \int_{D_L^c} dy\, \Gamma(t-s, x-y)\Bigg)\notag  \\
      &\qquad \leq \tilde C\, (1+ \Vert u_0 \Vert_\infty) \left(\sqrt{K_{t/2,x,L}} + K_{t,x,L}  \right) \notag \\
      &\qquad\qquad\qquad \times  \left[\exp\left(- \tfrac{(L - x)^2}{4 t} \right) + \exp\left(- \tfrac{(L + x)^2}{4 t} \right) \right] \notag \\
      &\qquad\leq  C_2\, (1+ \Vert u_0 \Vert_\infty)\, \sqrt{K_{t,x,L}} \, a_L(t/2, x),
 \label{s1-12}
\end{align}
because $K_{t,x,L} \leq \sqrt{K_{t,x,L}}$. From \eqref{bound0}--\eqref{s1-12}, we obtain \eqref{zero-one-two}.

Next, we study $A_3(t,x)$. Define
\begin{align*}
A_{3,1}(t,x)&= \int_0^t\int_{D_L} \Gamma(t-s,x-y)[\sigma(u(s,y))-\sigma(u_L(s,y))]\, W(ds,dy),\\
A_{3,2}(t,x)&= \int_0^t ds \int_{D_L} dy\,\Gamma(t-s,x-y)[b(u(s,y))-b(u_L(s,y))],
\end{align*}
so that $A_3(t,x) = A_{3,1}(t,x)+A_{3,2}(t,x)$. Applying Burkholder's inequality and using the Lipschitz property in ($\mathbf H_{\sigma, b}$), we see that
\beqn
\Vert  A_{3,1}(t,x)\Vert_{L^p(\Omega)}^2 \le C_{3,1}\left\Vert  \int_0^t ds \int_{D_L} dy\,\Gamma^2(t-s,x-y)(u(s,y)-u_L(s,y))^2\right\Vert_{L^{p/2}(\Omega)},
\eeqn
where $C_{3,1} = C_{3,1}(p, t, C_\tL)$. 
We apply Minkowski's inequality to bound from above the right-hand side 
and we obtain
\beqn
\label{a3.1}
\Vert  A_{3,1}(t,x)\Vert_{L^p(\Omega)}^2 \le C_{3,1}\int_0^t ds \int_{D_L} dy\,\Gamma^2(t-s,x-y) f_L^2(s,y).
\eeqn
Similarly, applying Minkowski's inequality to $A_{3,2}(t,x)$, we see that
\beq
\label{a3.2}
\Vert  A_{3,2}(t,x)\Vert_{L^p(\Omega)} \le C_\tL \int_0^t ds \int_{D_L} dy\,\Gamma(t-s,x-y) f_L(s,y).
\eeq
Writing $f_L(s,y) = 1\times f_L(s,y)$ and applying the Cauchy-Schwarz inequality to this integral, we obtain
\beq
\label{a3}
\Vert  A_{3}(t,x)\Vert_{L^p(\Omega)}^2\le \bar C_3\int_0^t ds \int_{D_L} dy\,\left[\Gamma^2(t-s,x-y)+ \Gamma(t-s,x-y)\right]f_L^2(s,y),
\eeq
where $\bar C_3 = \bar C_3(p, T, C_\tL)$.

Observe that
\beqn
   \Gamma(r, z) + \Gamma^2(r, z) 
    \leq \Gamma(r, z)\left(1 + \tfrac{1}{\sqrt{4 \pi r}}\right) , 
\eeqn
and for $r \in \, ]0, T]$ and $c = \sqrt{4 \pi T} + 1$,
\beqn
\Gamma(r, z) \left(1 + \tfrac{1}{\sqrt{4 \pi r}}\right) \leq  \tfrac{c}{\sqrt{4 \pi r}}\, \Gamma(r, z).
\eeqn
 Let
\beq\label{rd07_22e2}
J(r,z):= \tfrac{1}{\sqrt{4 \pi r}}\, \Gamma(r, z).
\eeq 
Notice that $\int_0^T dr \int_{\re} dz\, J(r,z) < \infty$.  
From the above computations, we deduce that 
\beq
\label{s1-13}
    \Vert A_3(t, x) \Vert_{L^p(\Omega)}^2 \leq C_3(p, T, C_\tL) \int_0^t ds \int_{D_L} dy \, J(t - s, x - y)  f_L^2(s, y).
\eeq
Then, using \eqref{s1-10}, \eqref{zero-one-two} and \eqref{s1-13}, and setting
\beqn
   b_L^2(t/2, x) = c\, (1 + \Vert u_0 \Vert_\infty)^2 \, K_{t,x,L}\, a_L^2(t/2, x),
\eeqn
we see that for $(t, x) \in \, ]0, T] \times D_L$,
\begin{align}
   f_L^2(t, x) &\leq  b_L^2(t/2, x) + C \int_0^t ds \int_{\R} dy \, J(t - s, x - y)\, 1_{D_L} (y)\, f_L^2(s, y).
\label{s1-14}
\end{align}
The constants $c$ and $C$ above depend on $p$, $T$ and $c_\tL$ and, in the case Nbc, also on $L_0 > 0$. Lemma \ref{a1-l3} in the Appendix (space-time Gronwall's Lemma, with $a(t, x) :=  b_L^2(t/2, x)$ and $f(t, x) := f_L^2(t, x)$)  yields, for $(t, x) \in \, ]0, T] \times D_L$,
\beq
\label{s1-15}
   f_L^2(t, x) \leq  b_L^2(t/2, x) +  [\cK \star b_L^2(\cdot/2, *)] (t, x),
\eeq
where the symbol ``$\star$" denotes space-time convolution and, for $(t, x) \in [0, T] \times \R$,
\beqn
   \cK(t, x) = \sum_{\ell = 1}^\infty C^\ell \, J^{\star \ell}(t, x)
\eeqn
and
\beq\label{rd07_21e1}
[\cK \star b_L^2(\cdot/2, *)] (t, x) = \int_0^t ds \int_{\R} dy \, \cK(t-s, x -y)\, b_L^2(s/2, y).
\eeq
We note that the assumptions of Lemma \ref{a1-l3} are satisfied because the functions $b_L^2(t/2, x)$, $f_L^2(t, x)$ and $ \cK(t, x)$ are bounded. 

The last part of the proof consist in showing that up to a positive multiplicative constant, $[\cK \star b_L^2(\cdot/2, *)] (t, x)$ is bounded by $(1 + \Vert u_0 \Vert_\infty)^2 \,a_L^2(t, x)$. 

A useful formula for  $\cK(t, x)$ is given in \cite[Proposition 2.2, Equation (2.8)]{chendalang2015-2}. More precisely, take $\lambda$ there equal to $\sqrt{C}$ and $\nu$ there equal to $4$, to see that
\beq
\label{boundfork} 
\cK(t,x) = \left(\tfrac{C}{4 \sqrt{\pi t}}  + \tfrac{C^2}{8} e^{C^2\, t/16} \Phi\left(C \, \sqrt{\tfrac{t}{8}} \, \right)\right) \Gamma(t,x),
\eeq
where $\Phi(\cdot)$ denotes the standard Normal probability distribution function. 

We bound $K_{t,x,L}$ by $1$ and omit for the moment the constant $c\, (1 + \Vert u_0 \Vert_\infty)^2$, to see that $(\cK \star b_L^2(\cdot/2, *)) (t, x)$ is bounded by the sum of two terms $I_\pm(t,x)$, where
\beqn
   I_\pm(t,x) = \int_0^t ds \int_{\R} dy \, \cK(t-s, x -y) \exp\left(-\tfrac{(L \pm y)^2}{2s} \right). 
    \eeqn
Use \eqref{boundfork} to see that
\begin{align*}
    I_\pm(t,x) &\leq \int_0^t ds \int_{\R} dy \, \cK(t-s, x -y)\, \sqrt{2 \pi s}\ \Gamma\left(s/2, L \pm y\right) \\
      & = \tilde C \int_0^t ds  \left(\tfrac{1}{\sqrt{t-s}} + e^{C^2\,(t-s)/16} \Phi\left(C\, \sqrt{\tfrac{t-s}{8}} \right) \right) \, \sqrt{s} \\
      &\qquad\qquad \times \int_{\R} dy \, \Gamma(t-s, x-y) \, \Gamma\left(s/2, L \pm y\right)\\
      &\leq  \tilde C' \int_0^t ds  \left(\tfrac{1}{\sqrt{t-s}} + 1\right) \, \sqrt{s}\ \Gamma(t - s/2, x \pm L),
   \end{align*}
 where we have used the semigroup property of the heat kernel. With the change of variables $r = t-s$, we obtain
\begin{align}\label{rd07_22e3}
 I_\pm(t,x) &\leq C_T \int_0^t dr  \left(\tfrac{1}{\sqrt{r}} +1 \right) \, \sqrt{t-r} \ \Gamma((t + r)/2, x \pm L) \\
    &= C_T \int_0^t dr  \left(\tfrac{1}{\sqrt{r}} +1 \right) \, \sqrt{t-r} \ \tfrac{1}{\sqrt{2 \pi (t + r)}} \exp\left(- \tfrac{(x \pm L)^2}{2(t + r)} \right) \notag \\
   & \leq \tilde C_T\,  e^{-\frac{(x\pm L)^2}{4t}}\int_0^t dr \left(\tfrac{1}{\sqrt r}+1\right) \sqrt{\tfrac{t-r}{t}} \notag \\
    &\le \tilde C_T'\,  e^{-\frac{(x\pm L)^2}{4t}}, 
\label{rd07_22e7} 
  \end{align} 
and consequently, for some constant $c_T$, for all $(t, x) \in \, ]0, T] \times D_L$,
\beqn
   [\cK \star b_L^2(\cdot/2, *)] (t, x)\le c_T\, (1+ \Vert  u_0\Vert_\infty)^2 \, a_L^2(t, x).
\eeqn
From the inequality \eqref{s1-15}, we conclude that $f_L(t,x)$ is bounded above by a constant times $(1 + \Vert u_0 \Vert_\infty)\, a_L(t,x)$. This completes the proof of  Theorem \ref{t-1}. 
\qed
\medskip


\section{Further special cases and optimality}
\label{rd08_25s1}

Using the same computations as in the proof of Theorem \ref{t-1} , we can revisit the particular examples considered in Section \ref{s1.0} and thereby gain additional insight into the optimality property mentioned in the Introduction.
We will see that in several special cases, the upper bound \eqref{s1-5} can be slightly improved, and a lower bound on the $L^2(\Omega)$-norm of the localization error can also be given under some assumptions on the coefficients $\sigma$ and $b$.


In the next Proposition, we study the deterministic case. This is an extension of Proposition \ref{rd05_04e1} to a nonlinear PDE.

\begin{prop}[The deterministic case $\sigma \equiv 0$]
\label{d08_25s1-p1}
The hypotheses are as in Theorem \ref{t-1}, and we assume that $\sigma \equiv 0$. Then the upper bound in \eqref{s1-5} can be improved to obtain
\begin{align}\nonumber
 & \Vert u(t, x) - u_L(t, x) \Vert_{L^p(\Omega)} \\
  &\qquad\qquad\qquad \leq c\, (1 + \Vert u_0 \Vert_\infty) 
   \left[\exp\left(-\tfrac{(L-x)^2}{4t} \right) + \exp\left(-\tfrac{(L+x)^2}{4t} \right) \right]
\label{rd07_22e4}
\end{align}
(notice that factors $4$ appear on the denominators instead of $8$).
\end{prop}
\begin{proof}
Going through the proof of Theorem \ref{t-1}, we see that in this case, $A_{1,2}(t,x)=0$
and the right-hand side of \eqref{zero-one-two} becomes
\beqn
    C_1\,  (1+ \Vert u_0 \Vert_\infty)\, K_{t, x, L} \,  a_L(t/2,x).
\eeqn
In addition, $A_{3,1} (t, x) \equiv 0$ and by \eqref{a3.2}, inequality \eqref{a3} is replaced by
\beqn
\Vert  A_{3}(t,x)\Vert_{L^p(\Omega)} \le \bar C_3\int_0^t ds \int_{D_L} dy\, \Gamma(t-s,x-y)\, f_L(s,y).
\eeqn
Formula \eqref{rd07_22e2} simplifies to
\beqn
    J(r, z) := \Gamma(r, z).
\eeqn
The inequality \eqref{s1-14} becomes
\begin{align*}
   f_L(t, x) &\leq  C_1\,  (1+ \Vert u_0 \Vert_\infty)\, K_{t, x, L} \,  a_L(t/2,x)  \\
   &\qquad\qquad + C \int_0^t ds \int_{D_L}  dy \, J(t - s, x - y)\,  f_L(s, y).
\end{align*}
The kernel $\cK(t, x)$ is now
\beq\label{rd05_05e4}
   \cK(t, x) = \Gamma(t, x) \sum_{\ell = 1}^\infty \frac{C^\ell}{(\ell - 1)!} =C \, e^C\,  \Gamma(t, x).
\eeq
Since, 
\begin{align*}
   [\cK \star a_L(\cdot/2, *)] (t, x) &\leq \int_0^t ds \int_\R dy \, C \, e^C\,  \Gamma(t-s, x-y)\, a_L(s/2, y) \\
   &\le c_T\,  a_L(t/2, x),
\end{align*}
if we bound $K_{t, x, L}$ by $1$, then, as a sharpening of \eqref{s1-15}, we obtain
\beqn
   f_L(t, x) \leq c \, (1+ \Vert u_0 \Vert_\infty)\, a_L(t/2, x),
\eeqn
and this establishes \eqref{rd07_22e4}. 
\end{proof}

 \begin{remark}
 By Proposition \ref{rd04_29e1}, the exponential terms in the upper bound in \eqref{rd07_22e4} are optimal, and the only possible improvement would be to obtain a pre-factor that depends on $p$, $t$, $L$ and $x$, for instance, a bound such as 
 \begin{align*}
  \Vert u(t, x) - u_L(t, x) \Vert_{L^p(\Omega)} 
  &\leq  c_{p}\, \sqrt{t}\, \left(1 + \Vert u_0 \Vert_\infty\right) \left(1 +  \left\vert\tfrac{L-x}{\sqrt{t}}\right\vert \right)^{-\frac34} \\
  & \qquad\qquad\times
   \left[\exp\left(-\tfrac{(L-x)^2}{4t} \right) + \exp\left(-\tfrac{(L+x)^2}{4t} \right) \right].
\end{align*}
 \end{remark}

The next Proposition considers additive noise. It shows that in this case, the exponential terms in 
Proposition \ref{rd04_29e1} can be recovered from the calculations in the proof of Theorem \ref{t-1}.

\begin{prop}[Additive noise $\sigma \equiv 1$]
\label{d08_25s1-p2}
The hypotheses are as in Theorem \ref{t-1}, and we assume that $\sigma \equiv 1$. Then the upper bound in \eqref{s1-5} can be improved to obtain
\eqref{rd07_22e4}.
\end{prop}
\begin{proof}
In the proof of Theorem \ref{t-1}, we see that, if $\sigma \equiv 1$ (or any constant) then the term $A_{3,1}(t,x)$ vanishes. The inequality
\eqref{a3} becomes 
\beq
\label{rd05_05e5}
\Vert  A_{3}(t,x)\Vert_{L^p(\Omega)}= \Vert  A_{3,2}(t,x)\Vert_{L^p(\Omega)} \le C_{\sigma, b} \int_0^t ds \int_{D_L} dy\, \Gamma(t-s,x-y)\, f_L(s,y),
\eeq
and \eqref{rd07_22e2} becomes 
 \beqn
     J(r, z) := \Gamma(r, z) .
 \eeqn
 The function $\cK(t,x)$ is now given by 
 \beqn
    \cK(t, x) = \sum_{\ell = 1}^\infty C^\ell \,   \Gamma^{\star \ell} (t, x) = \Gamma(t, x) \sum_{\ell = 1}^\infty \frac{C^\ell }{(\ell - 1)!} = C \, e^C\,  \Gamma(t, x).
\eeqn
We can now conclude in the same way as in Proposition \ref{d08_25s1-p1}.
\end{proof}

We end this section with the study of the case where there is no drift term.

\begin{prop}[Absence of drift $b \equiv 0$]
\label{d08_25s1-p3}
The hypotheses are as in Theorem \ref{t-1}, and we assume that $b \equiv 0$. Then the upper bound in \eqref{s1-5} can be improved to obtain
\eqref{rd07_22e4}.
\end{prop}
\begin{proof}
When $b \equiv 0$, the term $A_{3,2}(t,x)$ in the proof of Theorem \ref{t-1} vanishes. The inequality 
\eqref{a3} becomes 
\beq
\label{rd07_22e5}
\Vert  A_{3}(t,x)\Vert_{L^p(\Omega)}^2\le \bar C_3\int_0^t ds \int_{D_L} dy\, \Gamma^2(t-s,x-y)\, f_L^2(s,y),
\eeq
and \eqref{rd07_22e2} becomes 
 \beqn
     J(r, z) := \Gamma^2(r, z) = \tfrac{1}{\sqrt{8 \pi r}}\, \Gamma(r/2, z).
 \eeqn
 In order to obtain a formula for $\cK(t, x)$, we use  \cite[Proposition 2.2, Equation (2.8)]{chendalang2015-2} with $\lambda$ there equal to $\sqrt{C}$ and $\nu$ there equal to $2$, to see that
\beq
\label{rd07_22e6} 
   \cK(t,x) = \left(\tfrac{C}{\sqrt{8\pi t}}  + \tfrac{C^2}{4} e^{C^2\, t/8} \Phi\left(C \, \sqrt{\tfrac{t}{4}} \, \right)\right) \Gamma(t/2,x),
\eeq
The right-hand side of \eqref{rd07_22e3} becomes
\beqn
   C_T \int_0^t ds  \left(\tfrac{1}{\sqrt{s}} +1 \right) \, \sqrt{t-s} \ \Gamma(t/2, x \pm L),
\eeqn
which leads, instead of to \eqref{rd07_22e7},  to the inequality
\beqn 
    I_\pm(t, x) \leq C_T'\,  e^{-\frac{(x\pm L)^2}{2t}},
\eeqn
and this yields \eqref{rd07_22e4}.
\end{proof} 

Finally, in the case of linear multiplicative noise and the $L^2(\Omega)$-norm, we establish lower bound that matches those in Proposition \ref{rd04_29e1} (b) and (c). 

\begin{prop}[Linear multiplicative noise $\sigma(u) = u$, $b \equiv 0$]\label{rd06_12p1} 
The hypotheses are as in Theorem \ref{t-1}, and we assume that $\sigma(t, x, u) = \sigma(u) := u$ and $b \equiv 0$. Suppose that $m:= \inf_{y\in \re}u_0(y) \geq 0$. Then we have the following conclusions:

   (a) For $t > 0$ and $x \in \R$ fixed, there is $L_1 = L_1(t, x) > 0$ such that for all  $L > \max(L_1, \vert x \vert)$,  \begin{align}
 \label{mss05_22-2rd}
    \Vert u(t, x) - u_L(t, x) \Vert_{L^2(\Omega)} \geq m\, \tfrac{\sqrt{t}}{\sqrt{2 \pi}}\,  \left(1 + \tfrac{L-\vert x\vert}{\sqrt{t}} \right)^{-3/2} \exp\left(- \tfrac{(L-\vert x\vert)^2}{4t} \right) .
 \end{align}
 
   (b) Consider the case of vanishing Dirichlet boundary conditions Dbc. Then for $t > 0$ and $x \in \R$ fixed, there is $L_1 = L_1(t, x) > 0$ such that for all  $L > \max(L_1, \vert x \vert)$,  \begin{align}
 \label{mss05_22-3rd}
   & \Vert u(t, x) - u_L(t, x) \Vert_{L^2(\Omega)}\notag\\
   &\qquad\quad \ge  
 m \left[ \left(\sqrt{\tfrac{t}{\pi}}\, \tfrac{L-|x|}{2t+(L-|x|)^2}\right)^2 + \sqrt{\tfrac{t}{2\pi}}\,  \left(1 + \tfrac{L-\vert x\vert}{\sqrt{t}} \right)^{-3} \right]^{\frac12} 
  \exp\left(- \tfrac{(L-\vert x\vert)^2}{4t} \right) .
 \end{align}
\end{prop}

\begin{proof}
Using the notation in the proof of Theorem \ref{t-1}, we see that since $b \equiv 0$, and taking into account the domains of integration,
\begin{align*}
   \Vert u(t, x ) - u_L(t, x)\Vert_{L^2(\Omega)}^2 &= \left(A_0(t,x) + A_{1,1}(t, x)\right)^2 + E[A_2^2(t,x)] \\
   &\qquad + E[(A_{1, 2}(t, x) + A_3(t, x))^2]\\
    & \geq  \left(A_0(t,x) + A_{1,1}(t, x)\right)^2 + E[A_2^2(t,x)].
\end{align*}
Now $\left(A_0(t,x) + A_{1,1}(t, x)\right)^2$ matches the first square in \eqref{mss-5_22e2}, and
\begin{align*}
  E[A_2^2(t,x)] &= E\left[ \left(\int_0^t \int_{D_L^c} \Gamma(t-s, x-y) \, u(s, y)\, W(ds, dy)\right)^2\right] \\
     &= \int_0^t \int_{D_L^c} \Gamma^2(t-s, x-y) \, E[u^2(s, y)].
\end{align*}
By the Cauchy-Schwarz inequality, $E[u^2(s, y)] \geq (E[u(s, y)])^2$. Using \eqref{u} and the assumptions on $\sigma$ and $b$, we have $(E[u(s, y)])^2 = I_0^2(s, y) \geq m^2$ by \eqref{initial-cond} and the definition of $m$. We conclude that
 \begin{align*}
   \Vert u(t, x ) - u_L(t, x)\Vert_{L^2(\Omega)}^2 &\geq \left(A_0(t,x) + A_{1,1}(t, x)\right)^2  \\
   &\qquad + m^2  \int_0^t \int_{D_L^c} \Gamma^2(t-s, x-y).
\end{align*}

Conclusions (a) and (b) of this Proposition now follow the same lines as the proofs of parts (b) and (c) in Proposition \ref{rd04_29e1}.
\end{proof}

\begin{remark}[Uniformly positive $\vert \sigma \vert$, $b \equiv 0$]
If, instead of the assumption $\sigma(t, x, u) = \sigma(u) := u$, we require that for some $\sigma_0 > 0$ and for all $(t, x, u)$, $\vert \sigma(t, x, u) \vert \geq \sigma_0 > 0$, and maintain the requirement that $b \equiv 0$, then we also obtain the conclusions of Proposition \ref{rd06_12p1}, with $m$ in \eqref{mss05_22-2rd} replaced by $\sigma_0$, and with $m$ in \eqref{mss05_22-3rd} replaced by $\min(m, \sigma_0)$.
\end{remark}
 
\section{Comparing heat kernels on $[-L, L]$ and on $\R$}
\label{s3}

In order to estimate the difference $\vert \Gamma(t, x - y) - \Gamma_L(t; x, y) \vert$, we will make use of the following lemma concerning the Green's function with Neumann boundary conditions: it  provides a formula that will also be useful in the cases of Dirichlet and Mixed boundary conditions.

\begin{lemma}\label{rd08_16l1}
Let $\Gamma_L^N(t; x, y)$ be the Green's function defined in \eqref{3}, and let
$$
   H_L^N(t; x, y) = \Gamma(t, x-y) - \Gamma_L^N(t; x, y).
$$ 
Then $H_L^N(t; x, y) \leq 0$ and
\begin{align}\label{rd08_16e1}
    \int_{D_L} \vert H_L^N(t; x, y) \vert\, dy = \int_{D_L^c} \Gamma(t, x-y)\, dy.
\end{align}
\end{lemma}

\begin{proof}
The fact that $H_L^N(t; x, y)$ is negative follows from \eqref{l1-b0}. In order to establish \eqref{rd08_16e1}, we 
could now make use of \cite[Proposition 1.4.4]{d-ss-2024}, but we prefer to give a direct proof. Let $\Z^* = \Z \setminus \{0\}$. By the second formula in \eqref{3}, we see that
\begin{align}\nonumber
    \int_{D_L} \vert H_L^N(t; x, y) \vert\, dy &= \sum_{m \in \Z^*} \int_{-L}^L \Gamma(t, x - y + 4mL)\, dy \\
    &\qquad\qquad + \sum_{m \in \Z} \int_{-L}^L \Gamma(t, x+y+(4m+2)L)\, dy.
\label{rd08_16e2}
\end{align}
We use the change of variables $z = y - 4mL$ in the first integral, and $z = -y - (4m+2)L$ in the second integral, to see that this is equal to
\begin{align*}
  \sum_{m \in \Z^*} \int_{-L - 4mL}^{L - 4mL} \Gamma(t, x - z)\, dz + \sum_{m \in \Z} \int_{-3L - 4mL}^{-L - 4mL} \Gamma(t, x - z)\, dz.
\end{align*}
We split the first series into sums over $m \geq 1$ and over $m \leq -1$; in the sum over $m \geq 1$, we let $m = n+1$, and we regroup with the second series, to see that this is equal to
\begin{align*}
  &\sum_{n=0}^\infty \left(\int_{-5L - 4nL}^{-3L - 4nL} \Gamma(t, x - z)\, dz +  \int_{-3L - 4nL}^{-L - 4nL} \Gamma(t, x - z)\, dz \right) \\
  &\qquad + \sum_{m = -\infty}^{-1}  \left(\int_{-L - 4mL}^{L - 4mL} \Gamma(t, x - z)\, dz +  \int_{-3L - 4mL}^{-L - 4mL} \Gamma(t, x - z)\, dz \right)\\
  &= \sum_{n=0}^\infty \int_{-5L - 4nL}^{-L - 4nL} \Gamma(t, x - z)\, dz +  \sum_{m = -\infty}^{-1} \int_{-3L - 4mL}^{L - 4mL} \Gamma(t, x - z)\, dz,
 \end{align*}
 where, in each series, we have used the concatenation property of integrals. We observe that $-L - 4(n+1)L = -5L - 4nL$ and $L - 4(m+1)L = -3L - 4mL$, so we again use the concatenation property to see that this is equal to
 \begin{align*}
    \int_{-\infty}^{-L} \Gamma(t, x - z)\, dz + \int_L^\infty \Gamma(t, x - z)\, dz = \int_{D_L^c}  \Gamma(t, x - z)\, dz .
  \end{align*}
 This proves  Lemma \ref{rd08_16l1}.
\end{proof}

We now present estimates of the difference $\vert \Gamma(t, x - y) - \Gamma_L(t; x, y) \vert$ in the three cases that we are considering, beginning with the case of vanishing Dirichlet boundary conditions Dbc. 

\begin{lemma}
\label{s3-l1}
Let $\Gamma_L^D(t;x,y)$ be the Green's function defined in \eqref{2} and let $K_{t,x,L}$ be as defined in \eqref{rd07_18e2}.
Fix $T > 0$. For $t \in [0, T]$ and $x,y\in [-L,L]$, let
\beqn
   H_L^D(t; x, y) = \Gamma(t, x-y) - \Gamma_L^D(t; x, y).
\eeqn
Then $H_L^D(t; x, y) > 0$ and
\beq
\label{cd-1}
   \int_{-L}^L dy\, H_L^D(t; x, y) \leq  K_{t, x, L} \left[\exp\left(-\tfrac{(L-x)^2}{4t} \right) + \exp\left(-\tfrac{(L+x)^2}{4t} \right) \right],
\eeq
\beq
\label{cd-2}
   \int_0^t ds \int_{-L}^L dy \,  H_L^D(s; x, y)  \leq t\,K_{t, x, L}  \left[\exp\left(-\tfrac{(L-x)^2}{4t} \right) + \exp\left(-\tfrac{(L+x)^2}{4t} \right) \right]
\eeq
and
\beq
\label{cd-3}
   \int_0^t ds \int_{-L}^L dy \,  [H_L^D(s; x, y)]^2  \leq \sqrt{\tfrac{t}{\pi}}\, K_{t, x, L} \left[\exp\left(-\tfrac{(L-x)^2}{4t} \right) + \exp\left(-\tfrac{(L+x)^2}{4t} \right) \right]^2.
\eeq

\end{lemma}

\begin{proof}
The nonnegativity of $H_L^D$ follows from \eqref{l1-a1}. Let $(t, x) \in [0, T] \times [-L,L]$. By \eqref{2},
\beqn
 0 <  H_L^D(t; x, y) = - \sum_{m \in \Z^*}\Gamma(t,y-x+4mL) + \sum_{m \in \Z} \Gamma(t,x+y+(4m + 2)L),
\eeqn
therefore,
\begin{align*}
 0 <  \int_{-L}^L dy\, H_L^D(t; x, y) &\leq \sum_{m \in \Z^*}  \int_{-L}^L dy\, \Gamma(t,y-x+4mL) \\
    &\qquad\qquad + \sum_{m \in \Z}  \int_{-L}^L dy\, \Gamma(t,x+y+(4m + 2)L) .
\end{align*}
This is the same expression as in \eqref{rd08_16e2}, therefore, by Lemma \ref{rd08_16l1},   
\begin{align*}
   \int_{-L}^L dy\, H_L^D(t; x, y) &\leq \int_{D_L} dy\, \vert H_L^N(t; x, y)\vert = \int_{D_L^c} dy\, \Gamma(t; x - y).
\end{align*}
In particular, \eqref{cd-1} follows from \eqref{a1-l2-1}.

For the proof of \eqref{cd-2}, it suffices to notice that for $s \leq t$, $\exp(-z^2/(4s)) \leq \exp(-z^2/(4t))$, and to integrate \eqref{cd-1} from $0$ to $t$.

 In order to check \eqref{cd-3}, we first prove that if $t\in[0,T]$ and $x, y \in [-L, L]$, then
\beq
\label{cd-5}
    H_L^D(t; x, y) \leq \Gamma(t, x-L) + \Gamma(t, x+L).
\eeq
Indeed, from the second equality in \eqref{2}, we see that
\begin{align}
\label{cd-5-bis}
   H_L^D(t; x, y) 
   & = \Gamma(t,x-y)-\Gamma_L^D(t;x,y)\notag\\
   &= \Gamma(t,x+y+2L)\notag\\
   &\qquad -\sum_{m\in\Z^*}[\Gamma(t,x-y+4mL)-\Gamma(t,x+y+(4m+2)L)].
   \end{align}
   Applying the change of variables $n = -m-1$, we obtain
   \beqn
   \sum_{m= -\infty}^{-1} \Gamma(t,x+y+(4m+2)L) = \Gamma(t,x+y-2L) +\sum_{n=1}^\infty \Gamma(t,x+y-(4n+2)L),
   \eeqn
   and therefore, we deduce from \eqref{cd-5-bis} that
   \begin{align*}
   H_L^D(t; x, y) &=  \Gamma(t, x+ y + 2L) + \Gamma(t, x+ y - 2L) \\
    &\qquad - \sum_{m=1}^\infty [\Gamma(t, x - y + 4mL)- \Gamma(t, x+ y + (4m+ 2) L)]  \\
    &\qquad - \sum_{m=1}^{\infty} [\Gamma(t, x - y - 4mL) - \Gamma(t, x+ y - (4m + 2) L)].
\end{align*}
Both series are nonnegative (because for $m \geq 1$, $x + y + (4m+2) L \geq x - y + 4mL \geq 0$, and $x+y-(4m+2) L \leq x-y-4mL \leq 0$), so we remove them. 
Since $0 \leq x+L \le x+y+2L$ and $x+y-2L \leq x - L \leq 0$, we obtain \eqref{cd-5}.

Using \eqref{cd-5} and \eqref{cd-1}, it follows that for $s\in\, ]0,t]$,
\begin{align*}
\int_{-L}^L dy \,  [H_L^D(s; x, y)]^2 &\leq \sup_{z \in [-L, L]} H_L^D(s; x, z) \int_{-L}^L dy \,  H_L^D(s; x, y) \\
   &\leq  [\Gamma(s, x-L) + \Gamma(s, x+L)] \\
   &\qquad \times K_{s, x, L} \left[\exp\left(-\tfrac{(L-x)^2}{4s} \right) + \exp\left(-\tfrac{(L+x)^2}{4s} \right) \right] \\
   &= (4 \pi s)^{-\half} \, K_{s, x, L} \left[\exp\left(-\tfrac{(L-x)^2}{4s} \right) + \exp\left(-\tfrac{(L+x)^2}{4s} \right) \right]^2\\
   &\le  (4 \pi s)^{-\half}\, K_{t, x, L} \left[\exp\left(-\tfrac{(L-x)^2}{4t} \right) + \exp\left(-\tfrac{(L+x)^2}{4t} \right) \right]^2.
\end{align*}
Integrating from $s=0$ to $s=t$, we obtain \eqref{cd-3}.
\end{proof}


We now discuss the case of mixed boundary conditions Mbc.

\begin{lemma}
\label{s3-l3} 
Let $\Gamma_L^M(t;x,y)$ be the Green's function defined in \eqref{4}.Fix $T > 0$. Let $K_{t,x,L}$ be as defined in \eqref{rd07_18e2}. For $(t, x) \in [0, T] \times [-L,L]$, let
\beqn
   H_L^M(t; x, y) = \Gamma(t, x-y) - \Gamma_L^M(t; x, y).
\eeqn
Then 
\beq
\label{cn-10}
   \int_{-L}^L dy\, \vert H_L^M(t; x, y)\vert \leq K_{t, x, L} \left[\exp\left(-\tfrac{(L-x)^2}{4t} \right) + \exp\left(-\tfrac{(L+x)^2}{4t} \right) \right],
\eeq
\beq
\label{cn-11}
   \int_0^t ds \int_{-L}^L dy \,  \vert H_L^M(s; x, y)\vert   \leq t \, K_{t, x, L} \left[\exp\left(-\tfrac{(L-x)^2}{4t} \right) + \exp\left(-\tfrac{(L+x)^2}{4t} \right) \right]
\eeq
and
\beq
\label{cn-12}
   \int_0^t ds \int_{-L}^L dy \,  \left[H_L^M(s; x, y)\right]^2  \leq \sqrt{\tfrac{t}{\pi}}\, K_{t, x, L} \left[\exp\left(-\tfrac{(L-x)^2}{4t} \right) + \exp\left(-\tfrac{(L+x)^2}{4t} \right) \right]^2.
\eeq
\end{lemma}

\begin{proof}
Let $(t, x) \in [0, T] \times [-L,L]$. By \eqref{4},
\begin{align*}
 H_L^M(t; x, y) &= - \sum_{m \in \Z^*} (-1)^{m} \Gamma(t,y-x+4mL) \\
     &\qquad\qquad + \sum_{m \in \Z} (-1)^{m}  \Gamma(t,x+y+(4m + 2)L),
\end{align*}
therefore,
\begin{align*}
 \int_{-L}^L dy\, \vert H_L^M(t; x, y) \vert &\leq \sum_{m \in \Z^*}  \int_{-L}^L dy\, \Gamma(t,y-x+4mL) \\
    &\qquad\qquad + \sum_{m \in \Z}  \int_{-L}^L dy\, \Gamma(t,x+y+(4m + 2)L) .
\end{align*}
This is the same expression as in \eqref{rd08_16e2}, therefore, by Lemma \ref{rd08_16l1},   
\begin{align*}
   \int_{-L}^L dy\, \vert H_L^M(t; x, y)\vert &\leq \int_{D_L} dy\, \vert H_L^N(t; x, y)\vert = \int_{D_L^c} dy\, \Gamma(t; x - y).
\end{align*}
In particular, \eqref{cn-10} follows from \eqref{a1-l2-1}.

For the proof of \eqref{cn-11}, it suffices to notice that for $s \leq t$, $\exp(-z^2/(4s)) \leq \exp(-z^2/(4t))$, and to integrate \eqref{cn-10} from $0$ to $t$.

Turning to \eqref{cn-12}, we first check that if $t\in[0,T]$ and $x,y\in [-L,L]$, then $H_L^M(t;x,y)$ can be bounded above and below as follows:
\beq
\label{cn-14}
    - \Gamma(t, x+y - 2L) \leq H_L^M(t;x,y)\le \Gamma(t,x-L) + \Gamma(t,x+L).
\eeq

Indeed, the first inequality in \eqref{cn-14} follows from \eqref{r1-1-bis} and \eqref{l1-a1}, since 
 \begin{align*}
      \Gamma_L^M(t;x,y) &= \Gamma_{2L}^D(t; x-L, y-L) + \Gamma_{2L}^D(t; x-L, L-y) \\
      & \leq \Gamma(t, x-y) + \Gamma(t; x + y - 2L).
 \end{align*}
For the second inequality in \eqref{cn-14}, using \eqref{r1-1-bis}, and since $\Gamma_{2L}^D(t;x-L,L-y)\ge 0$, we obtain
\begin{align}
\label{cn-15}
H_L^M(t;x,y)&= \Gamma(t,x-y) - \Gamma_{L}^M(t;x,y)\notag\\
&\le \Gamma(t,x-y) - \Gamma_{2L}^D(t;x-L,y-L)\notag\\
&=\Gamma(t,x-y) \notag\\
&\qquad- \sum_{m\in\Z}[\Gamma(t,x-y+8mL)-\Gamma(t,x+y+2(4m+1)L)\notag]\\
&= \Gamma(t,x+y+ 2L)\notag\\
&\qquad -\sum_{m\in\Z^*}[\Gamma(t,x-y+8mL)- \Gamma(t,x+y+2(4m+1)L)]\notag†\\
&=: \Gamma(t,x+y+ 2L) - S_L(t,x,y).
\end{align}
Splitting the series $S_L(t,x,y)$ into the sum of two terms corresponding to positive and negative values of $m$, we see that
\begin{align*}
S_L(t,x,y) 
& = \sum_{m=1}^\infty [\Gamma(t,x-y+8mL)- \Gamma(t,x+y+2(4m+1)L)]\\
&\qquad + \sum_{m=1}^\infty[\Gamma(t,8mL+y-x)-\Gamma(t, x+y-2(4m-1)L)],
\end{align*}
where we have changed $m$ into $-m$ in the series corresponding to negative values of $m$. 
With the change of variables $k=m-1$, we have
\begin{align*}
 \sum_{m=1}^\infty \Gamma(t, x+y-2(4m-1)L) &= \sum_{k=0}^\infty \Gamma(t,x+y-2(4k+3)L) \\
 &= \Gamma(t,x+y-6L) + \sum_{m=1}^\infty\Gamma(t,2(4m+3)L-x-y),
 \end{align*}
 and consequently,
 \begin{align*}
 S_L(t,x,y)&= - \Gamma(t,x+y-6L)\\
 &\qquad +\sum_{m=1}^\infty[\Gamma(t,x-y+8mL)- \Gamma(t,x+y+2(4m+1)L)]\\
 &\qquad +\sum_{m=1}^\infty[\Gamma(t,8mL+y-x) - \Gamma(t,2(4m+3)L-x-y)].
 \end{align*}
 Because for $x,y\in [-L,L]$ and $m\ge 1$, $0<x-y+8mL\le x+y+2(4m+1)L$ and $0<8mL+y-x \le 2(4m+3)L-x-y$, we see that both series are nonnegative. Thus, $S_L(t,x,y)\ge  - \Gamma(t,x+y-6L)$ and from \eqref{cn-15} we deduce that
 \beqn
 H_L^M(t;x,y) \le\Gamma(t,x+y+2L)+ \Gamma(t,x+y-6L).
 \eeqn
 Finally, since for any $x,y\in[-L,L]$, $0\le x+L\le x+y+2L$ and $x+y-6L \le x-L\le 0$, we obtain the second inequality in \eqref{cn-14}. 
 
 By \eqref{cn-14},
 \begin{align}\nonumber
      \vert H_L^M(t; x, y)\vert &\leq \max\left(\Gamma(t, x+y - 2L),  \Gamma(t,x-L) + \Gamma(t,x+L) \right) \\
         &= \Gamma(t,x-L) + \Gamma(t,x+L), \label{rd08_17e1}
 \end{align}
 because $x + y - 2L \leq x-L \leq 0$, for all $x, y \in [-L, L]$.

 In order to prove \eqref{cn-12}, we observe the similarity between \eqref{rd08_17e1} and \eqref{cd-5}, and the similarity between \eqref{cn-11} and \eqref{cd-2}, and then we argue in the same way as for the proof of \eqref{cd-3} in Lemma \ref{s3-l1}.
 \end{proof}


We end this section with the study of the case of vanishing Neumann boundary conditions Nbc.

\begin{lemma}
\label{s3-l2} Let $\Gamma_L^N(t;x,y)$ be the Green's function defined in \eqref{3}.
Fix $T > 0$ and $L_0 > 0$. Let $K_{t,x,L}$ be as defined in \eqref{rd07_18e2}.  For $L > L_0$,  $t \in [0, T]$ and $x,y\in [-L,L]$, let
\beqn
   H_L^N(t; x, y) = \Gamma(t, x-y) - \Gamma_L^N(t; x, y).
\eeqn
Then $H_L^N(t; x, y) \leq 0$ and we have the following inequalities:
\beq
\label{cn-1}
   \int_{-L}^L dy\,\left\vert H_L^N(t; x, y)\right\vert \leq K_{t,x,L} \left[\exp\left(-\tfrac{(L-x)^2}{4t} \right) + \exp\left(-\tfrac{(L+x)^2}{4t} \right) \right],
\eeq
\begin{align}\nonumber
  & \int_0^t ds \int_{-L}^L dy \, \left\vert H_L^N(s; x, y)\right\vert \\
  &\qquad\qquad\qquad \leq t\, K_{t,x,L} \left[\exp\left(-\tfrac{(L-x)^2}{4t} \right) + \exp\left(-\tfrac{(L+x)^2}{4t} \right) \right],
\label{cn-2}
\end{align}
and
\begin{align}\nonumber
  & \int_0^t ds \int_{-L}^L dy \,  \left\vert H_L^N(s; x, y)\right\vert^2  \\
   &\qquad \leq \sqrt{\tfrac{t}{\pi}}\, 2 \vartheta\left(\tfrac{4L_0^2}{t}\right) K_{t,x,L}  \left[\exp\left(-\tfrac{(L-x)^2}{4t} \right) + \exp\left(-\tfrac{(L+x)^2}{4t} \right) \right]^2,
\label{cn-3}
\end{align}
where $\vartheta(a) = \sum_{m=0}^\infty e^{-am^2}$, $a>0$,  is the  {\em Theta function} (\cite[Chapter 5, Section 3.1]{stein-shak-2003}).
\end{lemma}

\begin{proof}
The non positivity of $H_L^N$ follows from \eqref{l1-b0}. The inequality \eqref{cn-1} is a direct consequence of Lemma \ref{rd08_16l1} and \eqref{a1-l2-1}.

For the proof of \eqref{cn-2}, it suffices to notice that for $s \leq t$, $\exp(-z^2/(4s)) \leq \exp(-z^2/(4t))$, and to integrate \eqref{cn-1} from $0$ to $t$.

 For the proof of \eqref{cn-3}, we use a method that could be applied to the other two types of boundary conditions, but would require taking $L_0 > 0$ there. We  first check that
\beq
\label{cn-4}
\left\vert H_L^N(t; x, y)\right\vert\le \tfrac{\vartheta(4L^2 /t)}{\sqrt{\pi t}} \left[\exp\left(-\tfrac{(L-x)^2}{4t}\right)  + \exp\left(-\tfrac{(L+x)^2}{4t}\right)\right].
\eeq 
Indeed, using \eqref{3} and rearranging terms, we see that
\begin{align}
\left\vert H_L^N(t; x, y)\right\vert &= \Gamma_L^N(t;x,y) - \Gamma(t,x-y)\notag\\
& = \sum_{m\in\Z^*} \Gamma(t,x-y+4mL) + \sum_{m\in\Z} \Gamma(t,x+y+(4m+2)L)\notag)\\
&= T_1 + T_2 + T_3 + T_4,
\label{cn-4-bis}
\end{align}
where
\begin{align*}
 T_1& =\sum_{m=1}^\infty\Gamma(t,x-y+4mL), \qquad\qquad T_2 = \sum_{m=-\infty}^{-1}\Gamma(t,x-y+4mL)\notag\\
 T_3&= \sum_{m=0}^\infty \Gamma(t,x+y+(4m+2)L), \qquad T_4 = \sum_{m=-\infty}^{-1} \Gamma(t,x+y+(4m+2)L).
\end{align*}
Let us first focus on the two series with $m\ge 0$. For $x,y\in [-L,L]$ and $m \geq 1$, we have 
$x-y+4mL \ge x+L+(4m-2)L\ge 0$, therefore
\begin{align*}
(x-y+4mL)^2 \geq (x+L+(4m-2)L)^2 \ge(x+L)^2 + (4m-2)^2L^2.
\end{align*}
We deduce that
\begin{align}
\label{cn-5}
T_1& 
   = \tfrac{1}{\sqrt{4\pi t}}\sum_{m=1}^\infty \exp\left(-\tfrac{(x-y+4mL)^2}{4t}\right)\notag\\
&\le \tfrac{1}{\sqrt{4\pi t}} \exp\left(-\tfrac{(x+L)^2}{4t}\right)\sum_{m=1}^\infty \exp\left(-\tfrac{(2m-1)^2L^2}{t}\right)\notag\\
&\le \tfrac{\vartheta(4L^2/t)}{\sqrt{4\pi t}} \exp\left(-\tfrac{(x+L)^2}{4t}\right).
\end{align}
Similarly, $x+y+(4m+2)L\ge x+L+4mL\ge 0$, therefore
\begin{align*}
     (x+L+4mL)^2 \ge (x+L)^2 + (4mL)^2
\end{align*}
and
\begin{align}
T_3& 
    = \tfrac{1}{\sqrt{4\pi t}}\sum_{m=0}^\infty \exp\left(-\tfrac{(x+y+(4m+2)L)^2}{4t}\right)\notag\\
&\le \tfrac{1}{\sqrt{4\pi t}} \exp\left(-\tfrac{(x+L)^2}{4t}\right)\sum_{m=0}^\infty \exp\left(-\tfrac{4m^2L^2}{t}\right)\notag\\
&\le  \tfrac{\vartheta(4L^2 /t)}{\sqrt{4\pi t}}\exp\left(-\tfrac{(x+L)^2}{4t}\right).
\label{cn-6}
\end{align}

Next, we consider the case $m\le -1$ and notice that for $x,y\in [-L,L]$ and $m \leq -1$, we have $x-y+4mL \le x-L+2(2m+1)L\le 0$, therefore
\begin{align*}
  (x-y+4mL)^2 \ge (x-L)^2 + 4(2m+1)^2L^2,
 \end{align*} 
 We deduce that
\begin{align}
T_2& 
   = \tfrac{1}{\sqrt{4\pi t}}\sum_{m=-\infty}^{-1}\exp\left(-\tfrac{(x-y+4mL)^2}{4t}\right)\notag\\
&\le \tfrac{1}{\sqrt{4\pi t}} \exp\left(-\tfrac{(x-L)^2}{4t}\right)\sum_{m=-\infty}^{-1}\exp\left(-\tfrac{(2m+1)^2L^2}{t}\right)\notag\\
&\le \tfrac{1}{\sqrt{4\pi t}} \exp\left(-\tfrac{(x-L)^2}{4t}\right)\sum_{m=0}^{\infty}\exp\left(-\tfrac{4m^2L^2}{t}\right)\notag\\
&= \tfrac{\vartheta(4L^2 /t)}{\sqrt{4\pi t}} \exp\left(-\tfrac{(x-L)^2}{4t}\right).
\label{cn-7}
\end{align}
Similarly, $x+y+(4m+2)L\le x-L+4(m+1)L\le 0$ for $m \leq -1$, therefore
\begin{align*}
    (x+y+(4m+2)L)^2  \ge  (x-L)^2 + 16(m+1)^2L^2
\end{align*}
and
\begin{align}
\label{cn-8}
T_4& 
     = \tfrac{1}{\sqrt{4\pi t}}\sum_{m=-\infty}^{-1}\exp\left(-\tfrac{(x+y+(4m+2)L)^2}{4t}\right)\notag\\
&\le \tfrac{1}{\sqrt{4\pi t}} \exp\left(-\tfrac{(x-L)^2}{4t}\right)\sum_{m=-\infty}^{-1}\exp\left(-\tfrac{4(m+1)^2L^2}{t}\right)\notag\\
&=\tfrac{\vartheta(4L^2 /t)}{\sqrt{4\pi t}} \exp\left(-\tfrac{(x-L)^2}{4t}\right).
\end{align}
From \eqref{cn-4-bis}-\eqref{cn-8}, we obtain \eqref{cn-4}.

Applying \eqref{cn-1} and \eqref{cn-4} we see that, for $s\in[0,T]$, 
\begin{align*}
&\int_{-L}^L dy\, \left\vert H_L^N(s;x,y)\right\vert^2\le  \sup_{z\in[-L,L]} \left\vert H_L^N(s;x,y)\right\vert\int_{-L}^L dy\, \left\vert H_L^N(s;x,y)\right\vert\\
&\qquad\le  \tfrac{\vartheta(4L^2 /s)}{\sqrt{\pi s}}\, K_{s, x, L}  \left[\exp\left(-\tfrac{(L-x)^2}{4s}\right)  + \exp\left(-\tfrac{(L+x)^2}{4s}\right)\right]^2\\
&\qquad\le \tfrac{1}{\sqrt s}  \tfrac{\vartheta(4L_0^2 /t)}{\sqrt{\pi}} \, K_{t, x, L}  \left[\exp\left(-\tfrac{(L-x)^2}{4t}\right)  + \exp\left(-\tfrac{(L+x)^2}{4t}\right)\right]^2,
\end{align*}
because the Theta function is nonincreasing and the exponentials above are increasing functions of the variable $s$. Integrating over the interval $[0,t]$, we obtain \eqref{cn-3}.
\end{proof}

\section{Appendix}
\label{a1}
In this section, we gather some facts concerning the Green's functions relative to the three types of boundary conditions described in Section \ref{s1}, and we also include some Gronwall-type lemmas that were used in the proof of Theorem \ref{t-1}. 
\medskip

\noindent{\em Auxiliary results on Green's functions}
\smallskip

The generic notation $\Gamma_L(t;x,y)$, used in the formulation of the equation for the random field solution \eqref{ul}, is replaced here by $\Gamma_L^D(t;x,y)$, $\Gamma_L^M(t;x,y)$ and  $\Gamma_L^N(t;x,y)$, in order to specify the type of boundary condition under consideration: Dirichlet, Mixed  and Neumann, respectively.
We recall below their expressions and their representations in terms of the heat kernel $\Gamma(t,z)$.
Notice that if in \eqref{s1-0}, we replace the space parameter $D_L=[-L,L]$ by $[0, L]$ and denote by $G_{L}^D(t;x,y)$ (respectively, $G_L^M(t;x,y)$ and $G_L^N(t;x,y)$) the Green's function associated to the heat equation on $[0, L]$ with Dirichlet (respectively, Mixed, Neumann)  boundary conditions, by a change of variables we see that $\Gamma_L^D(t; x, y)= G_{2L}^D(t;x+L,y+L)$ (respectively, $\Gamma_L^M(t; x, y)= G_{2L}^M(t;x+L,y+L)$, $\Gamma_L^N(t; x, y)= G_{2L}^N(t;x+L,y+L)$). 
\medskip

\noindent{\em The Green's function on $D_L$ for vanishing Dirichlet boundary conditions Dbc}
\begin{align}
\label{2}
\Gamma_L^D(t; x, y) &= \tfrac{1}{L} \sum_{n=1}^\infty e^{-\frac{\pi^2}{4 L^2}n^2 t} \sin\left(\tfrac{n \pi}{2L}(x+L) \right) \, \sin\left(\tfrac{n \pi}{2L}(y+L) \right)\notag \\
   &= \sum_{m \in \Z} \left[\Gamma(t,x-y+4mL) - \Gamma(t,x+y+(4m + 2)L)\right].
\end{align}
\smallskip

\noindent{\em The Green's function on $D_L$ for vanishing Mixed boundary conditions Mbc}
\begin{align}
\label{4}
\Gamma_L^M(t; x, y) &=  \tfrac{1}{L} \sum_{n=0}^\infty e^{- \frac{\pi^2 (2n+1)^2 t}{16 L^2}} \sin\left(\tfrac{(2n+1) \pi}{4L} (x + L) \right) \sin\left(\tfrac{(2n+1) \pi}{4L} (y + L) \right)\notag\\
   &= \sum_{m \in \Z} (-1)^m [\Gamma(t, x-y+ 4mL) - \Gamma(t, x+y +(4m+2)L)].
\end{align}
\smallskip

\noindent{\em The Green's function on $D_L$ for vanishing Neumann boundary conditions Nbc}
\begin{align}
\label{3}
\Gamma_L^N(t; x, y) &= \tfrac{1}{2L}+\tfrac{1}{L} \sum_{n=1}^\infty e^{-\frac{\pi^2}{4 L^2}n^2 t} \cos\left(\tfrac{n \pi}{2L}(x+L) \right) \, \cos\left(\tfrac{n \pi}{2L}(y+L) \right)\notag \\
   &= \sum_{m \in \Z} \left[\Gamma(t,x-y+4mL) + \Gamma(t,x+y+(4m+2)L)\right].
\end{align}


Detailed proofs of \eqref{2} and \eqref{3} are given in basic books on PDEs. They are also in \cite[Section 1.4]{d-ss-2024}. 
For  mixed boundary conditions, we refer to \cite[Section A.4.3]{candil2022}; the first equality in \eqref{4} follows from the identity (A.68) and the second from (A.72) there. 

\begin{remark}
\label{r1}
For any $t\in[0,T]$ and  $x,y\in [-L,L]$,
the Green's functions $\Gamma_L^M(t;x,y)$ and $\Gamma_{L}^D(t;x,y)$ are related. 
Indeed, using the first or second expressions in \eqref{2} and \eqref{4}, we see that for $t\in[0,T]$ and $x,y\in[-L,L]$,
\begin{align}
\label{r1-1-bis}
\Gamma_L^M(t;x,y) 
& = \Gamma_{2L}^D(t;x-L,y-L) + \Gamma_{2L}^D(t;x-L, L-y).
\end{align}
\end{remark}

\begin{lemma}
\label{l1}
\begin{description}
\item{(a)} The Green's function $\Gamma_L(t;x,y)$ relative to any one of the boundary conditions Dbc, Mbc and Nbc has the {\em semigroup property}: 
For any $s,t>0$ and $x,z\in D_L$,
\beq
\label{l1-semigroup}
\int_{D_L} dy\, \Gamma_L(s;x,y)\, \Gamma_L(t;y,z) = \Gamma_L(s+t;x,z).
\eeq
\item{(b)} {\em Case Dbc}. For any $t>0$ and $x,y\in D_L$, the Green's function $\Gamma_L^D(t;x,y)$ given in \eqref{2} has the following properties:
\beq
\label{l1-a1}
0\le \Gamma_L^D(t;x,y) < \Gamma(t,x-y)\le \tfrac{1}{\sqrt{4\pi t}},
\eeq
hence
\beq
\label{l1-a4}
\sup_{x\in D_L}\int_{D_L} dy\, \Gamma_L^D (t;x,y)<\int_{\re} \Gamma(t,z)\ dz =1,
\eeq
\beq
\label{l1-a2}
\sup_{x\in D_L} \int_{D_L} dy\, \left[\Gamma_L^D(t;x,y)\right]^2 \le \tfrac{1}{\sqrt{8\pi t}},
\eeq
and therefore
\beq
\label{l1-a3}
\int_0^t dr \sup_{x\in D_L}\int_{D_L} dy\, \left[\Gamma_L^D(r;x,y)\right]^2 \le \left(\tfrac{t}{2\pi}\right)^\half.
\eeq

\item{(c)} {\em Case Mbc}. For any $t>0$ and $x,y\in D_L$, 
\beq
\label{l1-c1}
     0 \leq \Gamma_L^M(t;x,y) \leq \tfrac{1}{\sqrt{\pi t}}\, , 
\eeq 
\beq
\label{l1-c2}
\sup_{x\in D_L} \int_{D_L} dy\, \left[\Gamma_L^M(t;x,y)\right]^2 \le \tfrac{1}{\sqrt{2 \pi t}}\, ,
\eeq
and
\beq
\label{l1-c3}
\int_0^t dr \sup_{x\in D_L} \int_{D_L} dy\, \left[\Gamma_L^M(r;x,y)\right]^2 \le  \sqrt{\tfrac{2 t}{\pi}}\, .
\eeq
Furthermore,
\beq
\label{l1-c4}
\int_{D_L} dy\, \Gamma_L^M (t;x,y) < 1.
\eeq

\item{(d)} {\em Case Nbc}. For any $t> 0$ and $x, y\in D_L$, we have:
\beq
\label{l1-b0}
0<\Gamma(t,x-y)\le\Gamma_L^N(t;x,y)  \le \tfrac{1}{2 L} +\tfrac{1}{\sqrt{\pi t}}
\eeq
and
\beq
\label{l1-b1}
\int_{D_L} dy\, \Gamma_L^N (t;x,y)=1,
\eeq
therefore 
\beq
\label{l1-b0-bis}
\sup_{x\in D_L}\int_{D_L} dy\, \Gamma_L^N (t;x,y)=1.
\eeq

Furthermore, for any $t\in\,]0,T]$ and $x,y\in D_L$, we have
\beq
\label{l1-b5-bis}
\sup_{x\in D_L} \int_{D_L} dy\, \left[\Gamma_L^N(t;x,y)\right]^2\ \le \tfrac{1}{2 L} +\tfrac{1}{\sqrt{2 \pi t}}
\eeq
and
\beq
   \label{int_J_2-tri}
  \int_0^t dr\ \sup_{x\in D_L} \int_{D_L} dy\, \left[\Gamma_L^N(r;x,y)\right]^2 \le \tfrac{t}{2 L} +\sqrt{\tfrac{2 t}{\pi}}.
  \eeq
\item{(e)} For the three cases of boundary conditions Dbc, Mbc and Nbc, for any $t\in[0,T]$,
  \beq
  \label{int_J_2}
  \int_0^t dr\ \sup_{x\in D_L} \int_{D_L} \Gamma_L(r;x,y)\ dy\le t. 
  \eeq
\end{description}
\end{lemma}

\begin{proof} 
(a) For each instance of boundary conditions, \eqref{l1-semigroup} can be easily checked using the first equalities in \eqref{2}--\eqref{3} and computing the integral on the left-hand side of \eqref{l1-semigroup}.

 (b) The function $\Gamma_L^D(t;x,\cdot)$ is the transition function of a Brownian motion that starts at $x$ and is absorbed at the boundary points $-L$ and $L$ (see e.g. \cite[Section 2.8]{ks}). This implies the first two inequalities in \eqref{l1-a1}. The third inequality follows from \eqref{s1-3}. 
 The inequality in \eqref{l1-a4} is a consequence of \eqref{l1-a1}.
 The upper bounds \eqref{l1-a2} and \eqref{l1-a3} are obtained using the semigroup property \eqref{l1-semigroup}, then \eqref{l1-a1} and computing the integrals.
 
(c)  The nonnegativity of $\Gamma_L^M(t; x, y)$ follows from the fact that $\Gamma_L^M(t; x, \cdot)$ is the probability density function of a Brownian motion in $[-L, L]$ that starts at $x$, is absorbed at the boundary $-L$ and is reflected on the boundary $L$. For the second inequality in \eqref{l1-c1}, we use \eqref{r1-1-bis} and \eqref{l1-a1} to see that
\begin{align*}
   \Gamma_L^M(t;x,y)  \leq \Gamma(t,x-y) + \Gamma(t,x+y-2L) \le  \tfrac{1}{\sqrt{\pi t}}.
\end{align*}

 By the semigroup property,
 \beqn
 \int_{D_L} [\Gamma_L^M (r;x,y)]^2\, dy = \Gamma_L^M(2r;x,x).
 \eeqn
 Hence, applying \eqref{l1-c1}, we obtain
 \beqn
 \sup_{x\in D_L}  \int_{D_L} [\Gamma_L^M (r;x,y)]^2\, dy =  \sup_{x\in D_L}  \Gamma_L^M(2r;x,x) \leq \tfrac{1}{\sqrt{2 \pi r}}, 
 \eeqn
 which is \eqref{l1-c2}. Integrating over $[0,t]$ with respect to $r$ yields \eqref{l1-c3}.

 Using  \eqref{r1-1-bis}, we see that
 \begin{align*}
    \int_{D_L} \Gamma_L^M(t; x, y) \, dy &= \int_{-L}^L \Gamma_{2L}^D(t; x-L,y-L)\, dy + \int_{-L}^L \Gamma_{2L}^D(t; x-L,L-y)\, dy \\
    & =  \int_{-2L}^0 \Gamma_{2L}^D(t; x-L,z)\, dz + \int^{2L}_0 \Gamma_{2L}^D(t; x-L,z)\, dz \\
    &=  \int_{D_{2L}} \Gamma_{2L}^D(t; x-L,z)\, dz \\
    & < 1, 
 \end{align*}
 by \eqref{l1-a4}. This proves \eqref{l1-c4}.  

 (d) The first inequality in \eqref{l1-b0} follows from \eqref{s1-3}, and the second inequality in \eqref{l1-b0} follows from the second formula in \eqref{3}.
 From the first formula in \eqref{3}, we see that
 \beqn
 \Gamma_L^N(t;x,y) \le \tfrac{1}{2L}+\tfrac{1}{L} \sum_{n=1}^\infty e^{-\frac{\pi^2}{4L^2}n^2 t}.
 \eeqn
 Since
 \beqn 
 \sum_{n=1}^\infty  e^{-\frac{\pi^2}{4L^2}n^2 t} \le \int_0^\infty dx\, e^{-\frac{\pi^2}{4L^2} x^2 t} = \tfrac{L}{\sqrt{\pi t}},
 \eeqn
 we deduce that
 \beqn
  \Gamma_L^N(t;x,y) \le \tfrac{1}{2L} + \tfrac{1}{\sqrt{\pi t}}. 
  \eeqn
This yields the last inequality in \eqref{l1-b0}.

 The identity \eqref{l1-b1} is obtained by integrating  the first formula in \eqref{3}, and \eqref{l1-b0-bis} follows immediately from \eqref{l1-b1}. 
 
 
 By the semigroup property,
 \beqn
 \int_{D_L} [\Gamma_L^N (r;x,y)]^2\, dy = \Gamma_L^N(2r;x,x).
 \eeqn
 Hence, applying \eqref{l1-b0}, we obtain
 \beqn
 \sup_{x\in D_L}  \int_{D_L} [\Gamma_L^N (r;x,y)]^2\, dy =  \sup_{x\in D_L}  \Gamma_L^N(2r;x,x) \le \frac{1}{2L} + \tfrac{1}{\sqrt{2 \pi r}}, 
 \eeqn
 which is \eqref{l1-b5-bis}. Integrating over $[0,t]$ with respect to $r$ yields \eqref{int_J_2-tri}.
 
(e)  Apply \eqref{l1-a4} (respectively,  \eqref{l1-c4} and \eqref{l1-b1}) to obtain \eqref{int_J_2} for $\Gamma_L:=\Gamma_L^D$ (respectively, $\Gamma_L:=\Gamma_L^M$ and $\Gamma_L:=\Gamma_L^N$). In the case Nbc,  \eqref{int_J_2} is in fact an equality.
  \end{proof}

\noindent{\em Gaussian tail estimates}
\medskip

We include the following elementary lemma so that the reader can check the constants in Lemma \ref{a1-l2}. 

\begin{lemma} For $\sigma > 0$, let
\beqn
    f_{\sigma^2}(x) = \tfrac{1}{\sqrt{2 \pi \sigma^2}} \, \exp\left(- \tfrac{x^2}{2 \sigma^2} \right), \qquad x \in \R.
\eeqn
For $a > 0$, 
\beq\label{rd07_18e1}
    \int_a^\infty f_{\sigma^2}(x) \, dx \leq \half\, \min\left(1, \sqrt{\tfrac{2}{\pi}}\, \tfrac{\sigma}{a} \right) \exp\left(- \tfrac{a^2}{2\sigma^2} \right).
\eeq
\end{lemma}

\begin{proof}
We have
\begin{align}\nonumber
& \int_a^\infty e^{-x^2/(2 \sigma^2)}\, dx = \int_0^\infty e^{- (z+a)^2/(2 \sigma^2)}\, dz  \\
    &\qquad\qquad \leq e^{- a^2/(2 \sigma^2)} \int_0^\infty e^{- z^2/(2 \sigma^2)}\, dz 
     = \half\, \sqrt{2\pi \sigma^2} \, e^{- a^2/(2 \sigma^2)}.
\label{rdlemC.2.2-before-3}
 \end{align}
 
  To complete the proof of \eqref{rd07_18e1}, notice that
 \begin{align*}
 \int_a^\infty e^{-x^2/(2\sigma^2)} \, dx & 
    \le 
      \int_a^\infty dx\, \tfrac{x}{a}\, e^{-\frac{x^2}{2 \sigma^2}}
 = 
     \tfrac{\sigma^2}{a}\, e^{-\frac{a^2}{2 \sigma^2}}.
 \end{align*}
Together with \eqref{rdlemC.2.2-before-3}, this yields \eqref{rd07_18e1}.
\end{proof}


The following lemma was used in the proof of Theorem \ref{t-1}. 

\begin{lemma}
\label{a1-l2}
Let $\Gamma(t,x)$, $t\in\, ]0,\infty[$, $x\in\re$, be the heat kernel defined  in \eqref{s1-3}, and let $K_{t,x,L}$ be as defined in \eqref{rd07_18e2}. Then, for $x \in D_L$,
\beq
\label{a1-l2-1}
\int_{D_L^c} \Gamma (t,x-y)\, dy\le K_{t, x, L}  \left[\exp\left(-\tfrac{(L-x)^2}{4t}\right) + \exp\left(-\tfrac{(L+x)^2}{4t}\right)\right],
\eeq
and
\beq
\label{a1-l2-2}   
  \int_{D_L^c} \Gamma^2 (t,x-y)\, dy\le \tfrac{K_{t/2, x, L}}{\sqrt{8 \pi t}} \left[\exp\left(-\tfrac{(L-x)^2}{2t}\right) + \exp\left(-\tfrac{(L+x)^2}{2t}\right)\right].
  \eeq
  \end{lemma}
  
  \begin{proof}
  Notice that $\Gamma(t,x-y) = f_{2t}(x-y)$ and
  \begin{align*}
      \int_{D_L^c} \Gamma (t,x-y)\, dy &= \int_{- \infty}^{x-L} f_{2t}(y)\, dy + \int_{x+L}^{\infty} f_{2t}(y)\, dy.
  \end{align*}
  By \eqref{rd07_18e1}, this is bounded above by
  \begin{align*}
   &  \half \min\left(1, \sqrt{\tfrac{t}{\pi}}\, \tfrac{2}{L-x} \right) \exp\left(- \tfrac{(L-x)^2}{4 t} \right)
     +  \half \min\left(1, \sqrt{\tfrac{t}{\pi}}\, \tfrac{2}{L+x} \right) \exp\left(- \tfrac{(L+x)^2}{4 t} \right) \\
     &\qquad \leq \half \max\left[\min\left(1, \sqrt{\tfrac{t}{\pi}} \, \tfrac{2}{L-x}\right) , \min\left(1, \sqrt{\tfrac{t}{\pi}} \, \tfrac{2}{L+x} \right)\right] \\
     &\qquad \qquad  \qquad \times\left[\exp\left(-\tfrac{(L-x)^2}{4t}\right) + \exp\left(-\tfrac{(L+x)^2}{4t}\right)\right].
\end{align*}
  This proves \eqref{a1-l2-1}. 
    
  Notice that
\begin{align*}
      \int_{D_L^c} \Gamma^2 (t,x-y)\, dy &=  \tfrac{1}{\sqrt{8 \pi t}}  \int_{D_L^c} \Gamma (t/2,x-y)\, dy.
\end{align*}
By \eqref{a1-l2-1}, this is bounded above by
\begin{align*}
    \tfrac{1}{\sqrt{8 \pi t}} \, K_{t/2, x, L} \left[\exp\left(-\tfrac{(L-x)^2}{2t}\right) + \exp\left(-\tfrac{(L+x)^2}{2t}\right)\right].
\end{align*}
This proves \eqref{a1-l2-2}.
  \end{proof}
  \medskip
 
 The following lemma is used in the proof of Proposition \ref{rd04_29e1}.
 
 \begin{lemma}\label{rd04_29l1}
   Fix $t > 0$ and $x \in \R$. Then
\begin{align*}
   \lim_{L \to \infty} \frac{ \int_0^t dr\, \exp\left(- \tfrac{(L-x)^2}{2r} \right)}{~\tfrac{1}{(L-x)^2} \exp\left(- \tfrac{(L-x)^2}{2t} \right)~} = 2 t^2.
\end{align*}
 \end{lemma}
 
 \begin{proof}
 The fraction is equal to
 \begin{align*}
    (L-x)^2 \int_0^{t} dr\, \exp\left( - \tfrac{(L-x)^2}{2} \left(\tfrac{1}{r} - \tfrac{1}{t} \right) \right) .
 \end{align*}
 Use the change of variables $u = \tfrac{1}{r} - \tfrac{1}{t}$ to see that this becomes
  \begin{align*}
  (L-x)^2 \int_0^{+\infty} \tfrac{du}{(u+1/t)^2}\, \exp\left( - \tfrac{(L-x)^2}{2}\, u \right) .
  \end{align*}
  Now do the change of variables $w = \frac{(L-x)^2}{2}\, u$ to see that this is equal to
   \begin{align*}
  2 \int_0^{+\infty}  dw \,  \tfrac{1}{(2w (L-x)^{-2} + 1/t)^2}\, \exp(-w).
   \end{align*}
 Using the dominated convergence theorem, we conclude that the limit as $L \to \infty$ of this expression is
\begin{align*}
      2\, t^2\int_0^{+\infty}  dw \, \exp(-w) = 2\, t^2.
\end{align*}
This proves the lemma.
 \end{proof}

 \noindent{\em A Gronwall-type lemma}
 \smallskip
 
 In the proof of Lemma \ref{rd08_09l1}, we have used a version of Gronwall's Lemma that appears frequently in the literature on SPDEs, and, 
 in the proof of Theorem \ref{t-1}, we have used a space-time Gronwall-type Lemma. We present both lemmas below. 
 
  Fix $T>0$. For two functions $f, g: \, ]0,T] \to \re_+$, their {\em convolution} (in time) is the function $f \ast g:  \, ]0,T] \to \re_+ \cup \{ \infty \}$ defined for  $t \in \, ]0, T]$ by $(f \ast g)(t) = \int_0^t f(t-s)\, g(s)\, ds$. Let $J: \, ]0,T] \to \re_+$ be a nonnegative Borel function such that
 $\int_0^T J(s)\, ds < \infty$.
 Set 
 \beqn
 u(t) = \sum_{n=1}^\infty J^{\ast n}(t), \quad U(t)=\int_0^t u(s) ds,
 \eeqn
 where $J^{\ast n}$ denotes the $n$-fold convolution of $J$ with itself. According to \cite[Lemma C.1.2]{d-ss-2024}, $0\leq U(t) \leq U(T) < \infty$.

 The next statement is a particular case of \cite[Lemma C.1.3 (c)]{d-ss-2024}.
 
 \begin{lemma}
 \label{a1-l3-usual}
  {\em (A Gronwall-type lemma)}
 Let $f: \, ]0,T] \to \re_+$ be a nonnegative Borel function. Suppose that there are nonnegative constants $c_1$ and $c_2$ such that for all $t\in\, ]0,T]$,
 \beqn
 f(t) \le c_1 + \int_0^t [c_2 + f(s)]\, J(t-s)\,ds.
 \eeqn
 Then for all $t\in\, ]0,T]$,
 \beqn
 f(t) \le c_1 + (c_1  + c_2)\,U(t).
 \eeqn
 \end{lemma}
\medskip
 
 \noindent{\em A space-time Gronwall-type lemma}
 \medskip
 
   Fix $T>0$. For two functions $f, g: \, ]0,T] \times \R \to \re_+$, their {\em convolution} (in space-time) is the function $f \star g:  \, ]0,T] \times \R \to \re_+ \cup \{ \infty \}$ defined for $(t, x) \in \, ]0, T] \times \R$ by $(f \star g)(t,x) = \int_0^t ds \int_\R dy\,  f(t-s, x-y)\, g(s,y)$. 
 
Fix $T>0$ and let $J: \, ]0,T]\times \re\to \re_+$ be a nonnegative Borel function such that 
\beq\label{rd07_16e1}
\int_0^T ds \int_{\re} dy\, J(s,y)<\infty.
\eeq
Let $J^{\star n}$ be the $n$-fold space-time convolution of $J$ with itself\ and let
\beqn
   \cK(t,x) = \sum_{n=1}^\infty J^{\star n}(t,x).
\eeqn
In the case of interest to us in this paper, we will have
\beq\label{rd07_16e2}
   \int_0^T ds \int_\R dy \, \cK(s, y) < \infty,
\eeq
so we assume this property here.

 \begin{lemma}
 \label{a1-l3} 
 {\em (A space-time Gronwall-type lemma)}
 Let $f, a:\,]0,T]\times \re \to \re_+$ be nonnegative Borel functions satisfying
 \beq
 \label{a1-l3-1}
 f(t,x)\le a(t,x) + \int_0^t ds \int_{\re} dy\,J(t-s,x-y) \,f(s,y),
 \eeq
for all $(t,x)\in \, ]0,T]\times \re$. Then for all $(t,x)\in \, ]0,T]\times \re$ such that $\lim_{n\to\infty} [J^{\star n} \star f](t,x)=0$ and $[\cK\star a](t,x)<\infty$ (these conditions are satisfied in particular if $f$ and $a$ are bounded),
\beq
\label{a1-l3-2}
f(t,x) \le a(t,x) +[\cK \star a](t,x).
\eeq
 \end{lemma}  
 
 \begin{proof}
 The inequality \eqref{a1-l3-1} can be rewritten $ f(t, x)  \leq a(t, x) + [J \star f](t, x)$. Replacing the $f$ on the right-hand side by the right-hand side of \eqref{a1-l3-1}, we obtain
 \beqn
     f(t, x)  \leq a(t, x) + [J \star a](t, x) + [J \star J\star f](t, x),
 \eeqn
 and we see by induction that for $n \geq 2$,
  \beqn
     f(t, x)  \leq a(t, x) + \sum_{\ell = 1}^{n-1} [J^{\star \ell} \star a](t, x) + [J^{\star n} \star f](t, x).
 \eeqn
 As $n \to \infty$, the third term converges to $0$ by assumption, and because $[\cK\star a](t,x)<\infty$, the second term converges to $[\cK\star a](t,x)$ by Fubini's theorem.
 \end{proof}
\bigskip

\noindent{\sc Acknowledgement.} The authors thank  the referee for a careful reading of the manuscript and several useful suggestions. These led in particular to the pre-factors in \eqref{mss05_22-3} and the lower bounds in Proposition \ref{rd06_12p1}.

\end{document}

{\bf Author information:}
\vskip 16pt

{\scshape
David Candil}

Institut de Mathématiques

\'Ecole Polytechnique Fédérale de Lausanne (EPFL)

CH-1015 Lausanne

Switzerland

\smallskip

david.candil@alumni.epfl.ch

\vskip 16pt

{\scshape
Robert C.~Dalang}

Institut de Mathématiques

\'Ecole Polytechnique Fédérale de Lausanne (EPFL)

CH-1015 Lausanne

Switzerland
\smallskip 

robert.dalang@epfl.ch

\vskip 16pt

{\scshape
Marta Sanz-Sol\'e}

Facultat de Matemàtiques i Informàtica

Universitat de Barcelona

Gran Via de les Corts Catalanes 585

08007 Barcelona

Spain
\medskip

and
\medskip

Royal Academy of Sciences and Arts of Barcelona

La Rambla 115

08002 Barcelona

Spain
\smallskip 

marta.sanz@ub.edu


\begin{thebibliography}{AA}

\bibitem{davies} B. Davies, Heat kernels and spectral theory. Cambridge Tracts in Mathematics, \textbf{92}. Cambridge University Press, Cambridge, 1989.

 \bibitem{candil2022} D.~J.-M.~Candil, Localization errors of the stochastic heat equation. Thèse no \textbf{7742}. \'Ecole Polytechnique F\'ed\'erale de Lausanne, 2022.
 
 \bibitem{CCL} D.~Candil, L. Chen and C.Y. Lee, {\em Parabolic stochastic PDEs on bounded domains with
rough initial conditions: moment and correlation bounds},  Stoch. Partial Differ. Equ. Anal. Comput.
\textbf{12} (3) (2024), 1507–1573.
 
 \bibitem{chendalang2015-2} L. Chen and R.C. Dalang,  {\em Moments and growth indices for the nonlinear stochastic heat equation with rough initial conditions}, Ann. Probab. \textbf{43}~(6) (2015), 3006–3051.


\bibitem{d-ss-2024} R.C. Dalang and M. Sanz-Sol\'e,  Stochastic partial differential equations, space-time white noise and random fields. Monographs in Mathematics. Springer, Cham, 2026.

\bibitem{gasquet} C. Gasquet and P. Witomski, Fourier analysis and applications. Filtering, numerical computation, wavelets. Texts in Applied Mathematics, \textbf{30}. Springer-Verlag, New York, 1999.


\bibitem{G-G-2017} M. Gerencsér and I. Gyöngy, {\em Localization errors in solving stochastic partial differential equations in the whole space},
Math. Comp. \textbf{86}  (307) (2017),  2373–2397.

\bibitem{G-G-2018} M. Gerencsér and I. Gyöngy,  {\em A Feynman-Kac formula for stochastic Dirichlet problems},
Stochastic Process. Appl. \textbf{129} no. 3 (2019), 995–1012.

\bibitem{H-R-S-W-2013} N. Hilber, O.  Reichmann, C. Schwab and C. Winter, {\em Computational methods for quantitative finance: Finite element methods for derivative pricing}. Springer Finance, Springer, Heidelberg, 2013.

\bibitem{ks} I. Karatzas and S. Shreve, Brownian motion and stochastic calculus. Second edition. Graduate Texts in Mathematics, \textbf{113}. Springer-Verlag, New York, 1991.

\bibitem{L-L-1996} D. Lamberton and B. Lapeyre, Introduction to stochastic calculus applied to finance.
Chapman \& Hall, London, 1996.

\bibitem{OLBC} F.W.J. Olver, D.W. Lozier, R.F. Boisvert, C.W. Clark (eds.). NIST Handbook of Mathematical Functions. U.S. Department of Commerce Institute of Standards and Technology. Washington DC, 2010.

\bibitem{stein-shak-2003} E.M. Stein and R. Shakarchi, Fourier analysis: an introduction.  Princeton Lectures in Analysis, \textbf{1}. Princeton University Press, Princeton, NJ, 2003.

\bibitem{walsh} J.B. Walsh,  
An introduction to stochastic partial differential equations.
 In {\em \'{E}cole d'\'et\'e de probabilit\'es de {S}aint-{F}lour,
  {XIV}---1984}, {\em Lecture Notes in Math.}  \textbf{1180}, pp. 265--439.
  Springer, Berlin  (1986).

\end{thebibliography}
\end{document}